\documentclass[final]{siamltex}
\usepackage{amsmath,amssymb,amsxtra,comment,graphicx,psfrag}
\usepackage{bm,mathrsfs}
\usepackage{mathtools}
\usepackage{xspace}
\usepackage{color}
\usepackage{enumerate,enumitem}
\usepackage{hyperref}
\usepackage{pdfsync}
\usepackage{hhline}
\usepackage{fancyhdr}
\usepackage{epsfig,subfigure,epstopdf}
\usepackage{pgfplots}

\allowdisplaybreaks

\topmargin \dimexpr 0.65in -0.99in
\newtheorem{remark}{\bf Remark}[section]
\newtheorem{example}{\bf Example}[section]

\def\u{{\bm u}}
\def\v{{\bm v}}
\def\H{{\bm H}}

\def\R{{\mathbb R}}

\def\d{{\mathrm d}}

\title{A convergent linearized Lagrange finite element method for the magneto-hydrodynamic equations\\ in 2D nonsmooth and nonconvex domains 
} 

\author{
Buyang Li\thanks{Department of Applied Mathematics, The Hong Kong Polytechnic University, Hung Hom, Hong Kong. E-mail address: buyang.li@polyu.edu.hk}
\and Jilu Wang\thanks{Department of Scientific Computing, 
Florida State University, Tallahassee, FL 32306, USA. E-mail address: jwang13@fsu.edu}
\and Liwei Xu\thanks{School of Mathematical Sciences, University of Electronic Science and Technology of China, Sichuan 611731, China. E-mail address: xul@uestc.edu.cn}
}

\begin{document}

\maketitle

\setlength\abovedisplayskip{2pt}
\setlength\belowdisplayskip{2pt}

\begin{abstract}
A new fully discrete linearized $H^1$-conforming Lagrange finite element method is proposed for solving the two-dimensional magneto-hydrodynamics equations based on a magnetic potential formulation. The proposed method yields numerical solutions that converge in general domains that may be nonconvex, nonsmooth and multi-connected. The convergence of subsequences of the numerical solutions is proved only based on the regularity of the initial conditions and source terms, without extra assumptions on the regularity of the solution. Strong convergence in $L^2(0,T;{\bf L}^2(\Omega))$ was proved for the numerical solutions of both ${\bm u}$ and ${\bm H}$ without any mesh restriction. 
\end{abstract} 

\begin{keywords}
{MHD, $H^1$-conforming, finite element, nonsmooth, nonconvex, convergence}
\end{keywords}

\begin{AMS} 
65M12, 
65N30, 
85A30  
\end{AMS}


\section{Introduction}\label{Se:intr} 
This article is concerned with numerical approximation of the incompressible magneto-hydrodynamics (MHD) equations 
\begin{align} 
&\mu\partial_t {\bm H}
+\sigma^{-1}\nabla\times(\nabla\times{\bm H})
-\mu\nabla\times({\bm u}\times{\bm H})
=\sigma^{-1}\nabla\times J 
\label{PDE-H}\\
&\partial_t{\bm u} +{\bm u}\cdot\nabla{\bm u}-\nu\Delta{\bm u}
+\nabla p 
={\bm f} - \mu {\bm H}\times(\nabla\times{\bm H})  \label{PDE-u-1} \\
&\nabla\cdot{\bm u}=0 \label{PDE-u-2} 
\end{align} 
in a two-dimensional polygonal type domain $\Omega=\Omega_0\backslash (\cup_{j=1}^m\Omega_j)\subset\R^2$, where both $\Omega_0$ and $\Omega_j\subset\Omega_0$, $j=1,\dots,m$, are polygons (thus the domain is possibly nonconvex and multi-connected), and the following 2D notations for the curl, divergence, and gradient operators are used for a vector field ${\bm B}=(B_1,B_2)$ and a scalar field $\psi$:  
\begin{align*}
&\nabla\times {\bm B}
=\frac{\partial B_2}{\partial x_1}-\frac{\partial B_1}{\partial x_2}
\qquad
\nabla\cdot {\bm B}=\frac{\partial B_1}{\partial x_1}
+\frac{\partial B_2}{\partial x_2} \\
&\nabla\times \psi=\bigg(\frac{\partial \psi}{\partial x_2} \,,\,
-\frac{\partial\psi}{\partial x_1}\bigg) \quad
\nabla\psi=\bigg(\frac{\partial \psi}{\partial x_1} \,,\,
\frac{\partial \psi}{\partial x_2}\bigg) \\
&
{\bm B}\times\psi=(B_2\psi\,,\,-B_1\psi) . 
\end{align*} 
In the system \eqref{PDE-H}-\eqref{PDE-u-2}, ${\bm u}$ denotes the velocity field, ${\bm H}$ the magnetic field, $p$ the pressure, ${\bm f}$ and $J$ the given source terms, $\nu$ the viscosity of the fluid, $\sigma$ the magnetic Reynolds number, and $\mu=M^2\nu\sigma^{-1}$, where $M$ denotes the Hartman number. 

We consider \eqref{PDE-H}-\eqref{PDE-u-2} under the perfectly conducting and no-slip boundary conditions 
\begin{align}\label{BC}
\H\cdot{\bm n}=0,
\quad 
\nabla\times\H = J ,
\quad\mbox{and}\quad
\u=0 
\quad\mbox{on}\,\,\,\partial\Omega\times(0,T] ,
\end{align} 
and the initial conditions
\begin{align}\label{ini}
\H|_{t=0}=\H_0
\quad\mbox{and}\quad 
\u|_{t=0}=u_0
\quad\mbox{in}\,\,\,\Omega .
\end{align}
The given source terms and the initial data for ${\bm H}$ and ${\bm u}$ are assumed to satisfy 
\begin{align}\label{Reg-f-J}
\begin{aligned}
&J\in C([0,T];L^2(\Omega))
\quad {\bm f}\in C([0,T];{\bf L}^2(\Omega)) \\ 
&
\H_0,\u_0\in {\bf L}^2(\Omega)\quad \quad \ \ 
\nabla\cdot\H_0=\nabla\cdot\u_0=0 . 
\end{aligned}
\end{align}
{\color{blue}
Under the assumption $\nabla\cdot\H_0 = 0$ for the initial value $\H_0$, a solution of \eqref{PDE-H} will automatically satisfy $\nabla\cdot\H = 0$ for all time.} 

The MHD equations \eqref{PDE-H}-\eqref{PDE-u-2} describe the interaction between a magnetic field and a viscous incompressible conducting fluid flow. {\color{blue}The mathematical theory of existence and uniqueness of weak solutions for the initial-boundary value problem \eqref{PDE-H}-\eqref{ini} was established in \cite[Theorem 3.1]{Sermange-Temam-1983} for a smooth or convex domain $\Omega$, where the weak formulation of the equations \eqref{PDE-H}-\eqref{ini} does not involve the variable $p$; see \cite[Problem 2.1]{Sermange-Temam-1983}.} 
In particular, under the regularity assumption \eqref{Reg-f-J} for the source terms and initial data, the problem has a unique weak solution (for any given $T>0$)
$$
({\bm H},{\bm u})\in L^2(0,T;{\bf H}^1(\Omega))^2 
\cap  L^\infty(0,T;{\bf L}^2(\Omega))^2 \quad\mbox{such that $\nabla\cdot{\bm H}=\nabla\cdot{\bm u}=0$}.
$$ 

Numerical methods and analysis for the MHD equations have been done from many different point of views. 
For the stationary MHD equations, existence, uniqueness, and finite element approximations were studied in \cite{GMP-1991} for small data. To overcome the numerical instability caused by possibly small hydrodynamic diffusion, a stabilized finite element method (FEM) was introduced in \cite{Gerbeau-2000}. These articles are concerned with $H^1$-conforming FEMs for the magnetic field ${\bm H}$ and the error estimates are based on the ${\bf H}^1(\Omega)$-regularity of the magnetic field ${\bm H}$. 
Such regularity holds in convex or smooth domains. However, 
in nonconvex and nonsmooth domains, the solution of the magnetic field is generally in ${\bf H}({\rm curl},\Omega)$ instead of ${\bf H}^1(\Omega)$. 

In more general domains, possibly nonconvex and nonsmooth, a mixed FEM with curl-conforming N\'ed\'elec edge elements was proposed for solving the magnetic field in the stationary MHD equations in \cite{Schotzau-2004}, where an additional gradient term $\nabla q$ was added to the magnetic potential equation to enforce the divergence-free condition for the magnetic field in a weak sense; the $H^1$-conforming FEM was used for the velocity. An error estimate for this numerical method was proved under the regularity assumptions 
\begin{align}\label{Reg-stationary}
({\bm u},p)\in {\bf H}^{s+1}(\Omega)\times H^s(\Omega)
\quad\mbox{and}\quad
({\bm H},\nabla\times{\bm H})\in {\bf H}^{s}(\Omega)\times H^s(\Omega)\quad\mbox{for some}\,\,\, s>\frac12 
\end{align}
which hold when the source term of the stationary MHD equations is sufficiently small. The same method for the magnetic field was also used in \cite{GLSW-2010} for solving the MHD equations, where a divergence-conforming FEM was used for the velocity. 

In the case of low magnetic Reynolds numbers, the MHD model usually consists of a time-dependent Navier-Stokes equation and a stationary electric potential equation with given magnetic field; {\color{blue}see \cite{GMP-1991,Peterson-1988}.} In \cite{Layton-Tran-Trenchea-2014}, two implicit-explicit methods (of first and second order) decoupling the velocity from electric potential  were proposed and analyzed. In \cite{RLZ-2018}, the method proposed in \cite{Layton-Tran-Trenchea-2014} was further combined with the technique of \cite{Shen-1995} to decouple the pressure from velocity. For the model with low magnetic Reynolds numbers considered in \cite{Layton-Tran-Trenchea-2014,RLZ-2018}, the $H^1$-conforming FEMs were proved to be convergent. 

For time-dependent MHD equations, many different numerical methods have been developed and analyzed:

\begin{enumerate}[label={\rm(\arabic*)},ref=\arabic*,topsep=2pt,itemsep=0pt,partopsep=1pt,parsep=1ex,leftmargin=20pt]

\item[$\bullet$]
$H^1$-conforming FEMs for the velocity and magnetic fields were used in \cite{Armero-Simo-1996,Guermond-Minev-2003,He-2014,He-Zou-2018,Zhang-He-2015} with several different time discretization methods.   These methods work very well in smooth or convex domains. In this type of  domains, error estimates have been established for several of these $H^1$-conforming FEMs. Similarly as for the stationary problem, in a nonconvex and nonsmooth domain, the solution of the magnetic field of the time-dependent MHD equations is generally in ${\bf H}({\rm curl},\Omega)$ in the spatial direction instead of ${\bf H}^1(\Omega)$. In this case the ${\bf H}^1(\Omega)$ or ${\bf H}^2(\Omega)$ regularity assumptions used in the existing analyses for the $H^1$-conforming FEMs generally do not hold. 

\item[$\bullet$]
A Galerkin least square FEM was proposed for solving the augmented MHD equations by adding an additional gradient term $\nabla q$ to the equation of magnetic field \cite{Salah-Soulaimani-Habashi-2001}, where numerical simulations were shown to illustrate the performance of the numerical methods. Rigorous proof for the convergence of numerical solutions remains open.

\item[$\bullet$]
The magnetic potential formulation was used in \cite{Shadid-2010}, where the equivalent formulation of the two-dimensional MHD equations 
\begin{align}
&\mu \partial_t A
-\sigma^{-1}\Delta A 
-\mu {\bm u}\times(\nabla\times A)=\sigma^{-1} J   \label{PDE-A0} \\
&\partial_t{\bm u} +{\bm u}\cdot\nabla{\bm u}-\nu\Delta{\bm u}
+\nabla p 
={\bm f} + \mu \nabla\cdot\big[(\nabla\times A)\otimes (\nabla\times A) -\frac12|\nabla\times A|^2\mathbb{I}\big]
\label{PDE-u-10} \\ 
&\nabla\cdot{\bm u}=0 \label{PDE-u-20} 
\end{align}
was solved by a fully implicit $H^1$-conforming FEM, where $\mathbb{I}$ denotes the identity matrix.  
Numerical results were given without proof for the convergence of numerical solutions. 

\item[$\bullet$]
Divergence-free preserving methods for MHD and ideal MHD have been developed in many articles.  
In particular, the locally divergence-free subspace was used in the discontinuous Galerkin (DG)  methods for the MHD equations \cite{Li-Shu-2005}. 
The divergence-free subspace of the Brezzi-Douglas-Marini  (BDM) finite element space was used for the magnetic field in \cite{Li-Xu-Yakovlev-2011,Li-Xu-2012}. 
An additional bubble function was added to each element in \cite{Cai-Wu-Xin-2013} in order to have additional degree of freedoms to enforce the divergence-free condition. 
More recently, the equation of magnetic field was rewritten as a first-order system in \cite{Hu-Ma-Xu-2017}, i.e., 
$$
\begin{aligned}
&\partial_t{\bm H}=-\nabla\times {\bm E} \\
&{\bm E}  = -{\bm u} \times {\bm H}+ \sigma^{-1}\nabla\times{\bm H} .
\end{aligned}
$$
Then, the divergence-conforming and curl-conforming N\'ed\'elec edge elements were used for ${\bm H}$ and ${\bm E}$, respectively. In this method no constraint was enforced on the magnetic field, but the numerical solution automatically satisfied the divergence-free condition provided that the initial data for the magnetic field was divergence free. 
For these methods, numerical results were shown to illustrate the performance of the numerical methods, while convergence proofs remain open. 

\item[$\bullet$]
{\color{blue}
A finite difference method for the 2D incompressible MHD equations in rectangular domains was proposed in \cite{Liu-Wang-2001}, where the MAC and Yee's scheme were used for the fluid and magnetic equations, respectively. The convergence of this method was proved in \cite{Ridder-2016} by using discrete energy estimates in $L^\infty(0,T;H^1(\Omega))$ and $L^2(0,T;H^2(\Omega))$ for both velocity and magnetic fields. 
An energy and helicity preserving finite difference method was proposed in \cite{Liu-Wang-2004} for the 3D axisymmetric incompressible MHD equations. The convergence of this method was demonstrated numerically but remains open theoretically.   
An error estimate of the corresponding energy and helicity preserving method for a single fluid equation was proved in \cite{Liu-Wang-2006} based on some regularity assumption.  
}

\item[$\bullet$]
As far as we know, the only existing proofs for the (subsequence) convergence of numerical solutions in possibly nonsmooth domains were given in \cite{Prohl-2008} and \cite{Hiptmair-Li-Mao-Zheng-2018}. 
In \cite{Prohl-2008}, Prohl studied several fully discrete linearized FEMs (with different time discretization and decoupling methods) with curl-conforming N\'ed\'elec edge elements for the magnetic field, and proved the convergence of two numerical schemes to weak solutions under the mesh restrictions $\tau=O(h^4)$ and $\tau=O(h^3)$, respectively, where $\tau$ denotes the time-step size and $h$ the spacial mesh size. Without such mesh restrictions, the weak$^*$ convergence of numerical solutions in $L^\infty(0,T;{\bf L}^2(\Omega))$ was proved. In \cite{Hiptmair-Li-Mao-Zheng-2018}, Hiptmair, Li, Mao and Zheng discretized a magnetic potential formulation of the three-dimensional MHD equations:
\begin{align} 
&\sigma \partial_t {\bm A} +\sigma (\nabla\times {\bm A})\times {\bm u} + \mu^{-1} \nabla\times \nabla\times{\bm A} = 0  \\
&\partial_t{\bm u} +{\bm u}\cdot\nabla{\bm u}-\nu\Delta{\bm u}
+\nabla p 
={\bm f}  -\sigma(\partial_t{\bm A}+(\nabla\times {\bm A})\times{\bm u})\times \nabla\times {\bm A} \\
&\nabla\cdot{\bm u}=0  
\end{align} 
which is different from the magnetic potential formulation in \cite{Shadid-2010}. 
The curl-conforming N\'ed\'elec elements were used for the discretization of ${\bm A}$. For every fixed time-step size $\tau$, it was proved that a subsequence $({\bm u}_{\tau,h_k}^n,{\bm H}_{\tau,h_k}^n)$, $k=1,2,\dots$, of numerical solutions will converge to a semi-discrete solution $({\bm u}_{\tau}^n,{\bm H}_\tau^n)$, i.e., 
$$
\begin{aligned}
&\lim_{k\rightarrow\infty}
\|{\bm u}_{\tau,h_k}^n-{\bm u}_{\tau}^n\|_{L^4(\Omega)}
\rightarrow 0 &&\mbox{for}\,\,\, n=1,\dots,N=T/\tau \\
&
\,\, {\bm A}_{\tau,h_k}^n \rightharpoonup {\bm A}_\tau^n 
\quad\mbox{weakly in $L^{\frac32}(\Omega)$ as $k\rightarrow\infty$} 
&&\mbox{for}\,\,\, n=1,\dots,N=T/\tau \\[5pt]
&
\lim_{k\rightarrow\infty}
\|{\bm H}_{\tau,h_k}^n-{\bm H}_\tau^n\|_{L^2(\Omega)}
\rightarrow 0 
&&\mbox{for}\,\,\, n=1,\dots,N=T/\tau. 
\end{aligned}
$$
It was also shown that a subsequence $({\bm u}_{\tau_m},{\bm A}_{\tau_m})$, $m=1,2,\dots$, of the semi-discrete solutions will converge weakly in $L^2(0,T;{\bf H}^1_0(\Omega))\times W^{1,p}(0,T;{\bf H}({\rm curl}))$ for some $p>1$ to a weak solution of the MHD equations. 
Strong convergence in $L^2(0,T;{\bf L}^2(\Omega))$ of the numerical solutions of ${\bm u}$ and ${\bf H}$ remains open. 

\end{enumerate}

Overall, the existing $H^1$-conforming FEMs for the MHD equations were proved to be convergent only in convex or smooth domains. In more general domains the mixed FEM with curl-conforming N\'ed\'elec elements is more suitable for approximating the magnetic field directly, while existing proofs for the convergence of numerical solutions either require mesh restriction $\tau=O(h^3)$ or yield only weak- or weak$^*$-convergence of the numerical solutions. 

{\color{blue}The most important applications of the MHD model considered in the present paper occur in metallurgy and liquid-metal processing \cite{Asai-2012, Unger-Mond-Branover-1988}. In that context, consideration of general polygonal yet nonconvex and multi-connected domains is very relevant.} The aim of this article is to develop a new fully discrete linearized $H^1$-conforming Lagrange FEM for the two-dimensional MHD equations based on a magnetic potential formulation such that the numerical solutions would converge not only in convex and smooth domains but also in nonconvex and nonsmooth domains. 
{\color{blue}
In the magnetic potential formulation, the magnetic potential $A$ would naturally has $H^1$ regularity and therefore can be solved correctly by using $H^1$-conforming finite element methods. Then the magnetic field is obtained through taking partial derivatives, i.e., $\H=\nabla \times A 
+ \sum_{j=1}^m\beta_j \nabla\times\varphi_j$, where $\beta_j$ and $\varphi_j$ are some time-independent constants and functions, depending on the geometry of the computational domain. 
This approach is different from the singular basis methods  \cite{Assous-Ciarlet-Sonnendrucker-1998,Hazard-Lohrengel-2003}, weighted regularization method \cite{Costabel-Dauge-2002}, $L^2$ projected method  \cite{DJLT-2009} and the residual-based stabilized augmented formulation \cite{Badia-Codina-2012}. These methods focus on approximating the magnetic field ${\bm H}$ in the Maxwell equations directly by using $H^1$-conforming finite element methods.} 
Similarly as \cite{Hiptmair-Li-Mao-Zheng-2018,Prohl-2008}, the proof of convergence in this paper is only based on the regularity of the initial conditions and source terms, without any extra assumptions on the regularity of the solution. Strong convergence of subsequences in $L^2(0,T;{\bf L}^2(\Omega))$ as $(\tau_n, h_n)\rightarrow (0,0)$ is proved for the numerical solutions of ${\bm u}$ and ${\bm H}$ without mesh restrictions.

The rest of the paper is organized as follows. In Section 2, we present an equivalent magnetic potential formulation of the the two-dimensional MHD equations \eqref{PDE-H}-\eqref{PDE-u-2}. In Section 3, we propose a fully discrete linearized $H^1$-conforming Lagrange finite element method for solving the problem, and present the main theoretical result about the convergence the numerical solutions. 
Rigorous proof of the main theoretical result is presented in Section 4. 
Numerical experiments are given in Section 5 to support the theoretical analyses.

\section{Equivalent formulation}
In section \ref{Sec:Formal} we first formally derive an equivalent formulation of the two-dimensional MHD equations \eqref{PDE-H}-\eqref{PDE-u-2} in terms of the magnetic potential. Then we define weak solutions of the problem in section \ref{Sec:weak}.  
It is easy to verify that a weak solution of the reformulated problem is also a weak solution of the original MHD equations (see Remark 2.1). 

\subsection{Formal derivation}\label{Sec:Formal} 
By taking the divergence of \eqref{PDE-H} we obtain 
$\mu\partial_t\nabla\cdot\H=0,$
which together with the divergence-free initial condition $\nabla\cdot{\bm H}_0=0$ in \eqref{Reg-f-J} implies 
\begin{align}\label{div-H=0}
\nabla\cdot\H=0.
\end{align} 
Let $m$ denote the number of holes of the domain $\Omega$, and let $\Gamma_j=\partial\Omega_j$ denote the boundary of the $j$th hole. Then the divergence-free vector field $\H$ can be decomposed as (cf. \cite{BCNS-2012}, with slightly different boundary conditions) 
\begin{align}\label{Hodge-D}
&\H=\nabla \times A 
+ \sum_{j=1}^m\beta_j \nabla\times\varphi_j,
\end{align}
where 
$\beta_j$ for $j=1,\dots,m$, are constants independent of the spatial variables, 
$A$ is the solution of 
\begin{align}\label{hodge-AB}
\left\{\begin{aligned}
&-\Delta A=\nabla\times\H &&\mbox{in}\,\,\,\Omega \\
&A=0 &&\mbox{on}\,\,\,\partial\Omega,
\end{aligned}\right.
\end{align}
and $\varphi_j$ is the solution of 
\begin{align}\label{hodge-varphi-j}
\left\{\begin{aligned}
&\Delta \varphi_j=0 &&\mbox{in}\,\,\,\Omega \\
&\varphi_j=1&&\mbox{on}\,\,\,\Gamma_j \\
&\varphi_j=0 &&\mbox{on}\,\,\,\partial\Omega\backslash\Gamma_j. 
\end{aligned}\right.
\end{align}

{\color{blue}Note that for the scalar-valued function $\phi$, there holds 
$$
\nabla\times\nabla\times\phi
=\nabla\times \bigg(\frac{\partial \phi}{\partial x_2} \,,\,
-\frac{\partial\phi}{\partial x_1}\bigg) 
=
\frac{\partial}{\partial x_1}\bigg(-\frac{\partial \phi}{\partial x_1}\bigg)
-\frac{\partial}{\partial x_2}\bigg(\frac{\partial \phi}{\partial x_2}\bigg)
=-\Delta\phi . 
$$
}
Therefore, integrating the time derivative of \eqref{Hodge-D} against $\nabla\times\varphi_i$ yields 
\begin{align*} 
&(\partial_t\H,\nabla\times\varphi_i) = \sum_{j=1}^m \frac{\d\beta_j}{\d t} (\nabla\times\varphi_j,\nabla\times\varphi_i)  ,
\end{align*}
where we have used \eqref{hodge-AB}-\eqref{hodge-varphi-j} and 
$$
(\nabla\times A,\nabla\times\varphi_i)
=-(A,\nabla\times(\nabla\times\varphi_i))
=(A,\Delta\varphi_i)
=0.
$$
Furthermore, integrating \eqref{PDE-H} against $\nabla\times\varphi_i$ and using integration by parts, and with the boundary condition \eqref{BC}, we obtain 
\begin{align*} 
(\partial_t\H,\nabla\times\varphi_i) =0 \quad\mbox{for $i=1,\dots,m$} .  
\end{align*}
Since the matrix $M_{ij}=(\nabla\times\varphi_j,\nabla\times\varphi_i) $ is positive definite, the two identities above imply 
\begin{align} 
\frac{\d\beta_j}{\d t} =0\quad\mbox{for}\ i=1,\dots,m .
\end{align} 
Therefore, $\beta_j$ for $j=1,\dots,m$, are constants independent of time.

Now we substitute $\H=\nabla\times A+ \sum_{j=1}^m\beta_j \nabla\times\varphi_j$ into \eqref{PDE-H}. This yields 
\begin{align*}
&\nabla\times \big(\mu \partial_t A
-\sigma^{-1}\Delta A -\sigma^{-1} J
-\mu{\bm u}\times(\nabla\times A+ \mbox{$\sum_{j=1}^m$}\beta_j \nabla\times\varphi_j)\big) 
={\bm 0}  
\end{align*}
which implies 
\begin{align}\label{eq-A}
&\mu \partial_t A -\sigma^{-1}\Delta A -\sigma^{-1} J
-\mu{\bm u}\times \big(\nabla\times A+ \mbox{$\sum_{j=1}^m$}\beta_j \nabla\times\varphi_j \big) = {\rm const}  . 
\end{align}
With the boundary condition \eqref{BC}, it is easy to derive that the constant on the right-hand side of the above equation equals zero.

Thus, instead of solving \eqref{PDE-H}-\eqref{PDE-u-1} directly, we propose to solve \eqref{hodge-varphi-j} and the following equations: 
\begin{align}
&\mu \partial_t A
-\sigma^{-1}\Delta A 
-\mu {\bm u}\times\big(\nabla\times A+ \mbox{$\sum_{j=1}^m$}\beta_j \nabla\times\varphi_j \big) =\sigma^{-1} J   \label{PDE-A} \\
&\partial_t{\bm u} +{\bm u}\cdot\nabla{\bm u}-\nu\Delta{\bm u}
+\nabla p 
={\bm f} + \mu\big(\nabla\times A+ \mbox{$\sum_{j=1}^m$}\beta_j \nabla\times\varphi_j \big) \times \Delta A 
\label{PDE-u-11} \\
&\nabla\cdot{\bm u}=0 , \label{PDE-u-22} 
\end{align}
where $\beta_j$, $j=1,\dots,m$, are determined by 
\begin{align}\label{betaj}
& \sum_{j=1}^m  \beta_j  (\nabla\times\varphi_j,\nabla\times\varphi_i) 
= (\H_0,\nabla\times\varphi_i) 
\end{align}
which can be obtained by integrating \eqref{Hodge-D} against $\nabla\times\varphi_i$ at time $t=0$. 

The boundary and initial conditions for \eqref{PDE-A}-\eqref{PDE-u-22} are given by 
\begin{align}\label{A-BC}
\u=0 , \quad 
A =0 
\quad\mbox{on}\,\,\,\partial\Omega\times(0,T] 
\end{align} 
and 
\begin{align}\label{A-ini}
\u|_{t=0}=u_0, 
\quad 
A|_{t=0}=A_0
\quad\mbox{in}\,\,\,\Omega ,
\end{align}
where $A_0$ is the solution of 
\begin{align}\label{Def-A0}
\left\{\begin{aligned}
&-\Delta A_0=\nabla\times\H_0 &&\mbox{in}\,\,\,\Omega \\
&A_0=0 &&\mbox{on}\,\,\,\partial\Omega.
\end{aligned}\right.
\end{align}
After solving of \eqref{hodge-varphi-j} and \eqref{PDE-A}-\eqref{A-ini}, we can obtain the magnetic field $\H=\nabla\times A+ \mbox{$\sum_{j=1}^m$}\beta_j \nabla\times\varphi_j$.

\medskip

\subsection{Weak solution}\label{Sec:weak}
For $k\ge 0 $ and $p\in[1,\infty]$, we denote by $W^{k,p}(\Omega)$ the conventional Sobolev space of functions defined on $\Omega$, with abbreviations $L^{p}(\Omega)=W^{0,p}(\Omega)$ and $H^{k}(\Omega)=W^{k,2}(\Omega)$. 
Let $W^{1,p}_0(\Omega)$ be the space of functions in $W^{1,p}(\Omega)$ with zero traces on the boundary $\partial\Omega$, and denote $H^1_0(\Omega)=W^{1,2}_0(\Omega)$. 
The corresponding vector-valued spaces are 
\begin{align*}
&{\bf L}^p(\Omega)=L^p(\Omega)\times L^p(\Omega)
&& {\bf W}^{k,p}(\Omega)=W^{k,p}(\Omega)\times W^{k,p}(\Omega)
\\
&{\bf W}^{1,p}_0(\Omega)=W^{1,p}_0(\Omega)\times W^{1,p}_0(\Omega) 
&& {\bf H}^1_0(\Omega)={\bf W}^{1,2}_0(\Omega) \\ 
&{\bf H}^1_{0,{\rm div}}(\Omega)=\{{\bm v}\in {\bf H}^1_0(\Omega):\nabla\cdot{\bm v}=0\} .
\end{align*}
The inner product of $L^2(\Omega)$ is denoted by $(\cdot,\cdot)$. 

\vspace{.1in}
A quadruple $(A,\u,(\varphi_j)_{j=1}^m,(\beta_j)_{j=1}^m)$ is called a weak solution of \eqref{hodge-varphi-j} and \eqref{PDE-A}-\eqref{Def-A0} if 
\begin{align}
&\!\!A\in L^\infty(0,T;H^1(\Omega)),\,\,\, \Delta A\in L^2(0,T;L^2(\Omega)),\,\,\,
\partial_tA\in L^s(0,T;L^2(\Omega)) 
\,\,\,\forall\, s\in(1,2) \label{RegA}\\
&\!\!\u\in L^\infty(0,T;{\bf L}^2(\Omega))\cap L^2(0,T;{\bf H}^1_{0,{\rm div}}(\Omega)), 
\,\,\,
\partial_t\u\in L^s(0,T;{\bf H}^1_{0,{\rm div}}(\Omega)') 
\,\,\,\forall\, s\in(1,2) \label{Regu}\\
&\!\!\varphi_j\in H^{\frac{3}{2}+\delta}(\Omega),\,\,\, j=1,\dots,m,
\,\,\,\mbox{are solutions of \eqref{hodge-varphi-j}, where $\delta\in(0,1)$ is a constant}  \\
&\!\!\beta_j\in\R,\,\,\, j=1,\dots,m,\,\,\,\mbox{satisfy}\,\,\,\eqref{betaj} \label{betaj2}
\end{align}
with $A|_{t=0}=A_0$ and $\u|_{t=0}=\u_0$, and the following equations hold for all test functions $a\in L^\infty(0,T;H_0^1(\Omega))$ and ${\bm v}\in L^\infty(0,T;{\bf H}^1_{0,{\rm div}}(\Omega))$, 
\begin{align}
&\int_0^T\Big[\big(\mu \partial_tA ,a\big) 
-\big(\sigma^{-1}\Delta A ,a \big)
- \big(\mu{\bm u}\times\big(\nabla\times A+ \mbox{$\sum_{j=1}^m$}\beta_j \nabla\times\varphi_j \big) ,a\big) \Big]\d t \nonumber \\
&=\int_0^T\big(\sigma^{-1} J ,a\big) \d t   \label{weak-A} \\[8pt] 
&\qquad\int_0^T \Big[\big(\partial_t\u,{\bm v}\big)
+\frac12({\bm u}\cdot\nabla{\bm u},{\bm v})
-\frac12({\bm u}\cdot\nabla {\bm v},{\bm u})
+(\nu\nabla{\bm u},\nabla {\bm v})
\nonumber \\
&=\int_0^T 
\Big[({\bm f},{\bm v})  + \big(\mu\big(\nabla\times A+ \mbox{$\sum_{j=1}^m$}\beta_j \nabla\times\varphi_j \big)\times \Delta A,{\bm v}\big)  \Big]\d t . 
\label{weak-u}
\end{align}

\noindent{\bf Remark 2.1.}
\begin{enumerate}[label={\rm(\arabic*)},ref=\arabic*,topsep=2pt,itemsep=0pt,partopsep=1pt,parsep=1ex,leftmargin=20pt]
\item
The pressure $p$ does not appear in the definition of weak solutions as we have restricted both $\u$ and the test function $\v$ to ${\bf H}^1_{0,{\rm div}}(\Omega)$. 

\item
If the domain is simply connected (without any holes) then $m=0$. In this case, a pair $(A,\u)$ is called a weak solution if \eqref{RegA}-\eqref{Regu} and \eqref{weak-A}-\eqref{weak-u} hold. 

\item
By substituting $a=\nabla\times{\bm b}$ into \eqref{weak-A}, we see that if $(A,\u,(\varphi_j)_{j=1}^m,(\beta_j)_{j=1}^m)$ is a weak solution of \eqref{hodge-varphi-j} and \eqref{PDE-A}-\eqref{Def-A0} then 
$({\bf H},\u)$ is a weak solution of \eqref{PDE-H}-\eqref{ini} in the sense that 
\begin{align}
&{\bm H}=\nabla\times A+ \mbox{$\sum_{j=1}^m$}\beta_j \nabla\times\varphi_j \in L^\infty(0,T;{\bf L}^2(\Omega))\cap L^2(0,T;{\bf H}(\Omega;{\rm curl})), \label{weak-H-reg} \\
&
\partial_t{\bm H}\in L^s(0,T;{\bf H}(\Omega;{\rm curl})') ,\quad \forall\, s\in(1,2),\\
&\mbox{$\H|_{t=0}=\H_0$ and $\u|_{t=0}=u_0$},
\end{align}
and the following equations hold for all test functions ${\bm b}\in L^\infty(0,T;{\bf H}(\Omega;{\rm curl}))$ and ${\bm v}\in L^\infty(0,T;{\bf H}^1_{0,{\rm div}}(\Omega))$, 
\begin{align}
&\int_0^T\Big[\big(\mu \partial_t\H,{\bm b}\big) 
+\big(\sigma^{-1}\nabla\times\H,\nabla\times{\bm b} \big)
- \big(\mu{\bm u}\times \H,\nabla\times{\bm b}\big) \Big]\d t
=\int_0^T\big(\sigma^{-1} J ,\nabla\times{\bm b}\big) \d t,   \label{weak-H} \\
&\int_0^T \Big[\big(\partial_t\u,{\bm v}\big)
+\frac12({\bm u}\cdot\nabla{\bm u},{\bm v})
-\frac12({\bm u}\cdot\nabla {\bm v},{\bm u})
+(\nu\nabla{\bm u},\nabla {\bm v})
\nonumber \\
&=\int_0^T 
\Big[({\bm f},{\bm v})  
-(\mu{\bm H} \times(\nabla\times {\bm H}),{\bm v})  \Big]\d t .   \label{weak-uu} 
\end{align}

\end{enumerate}

\section{Numerical method}\label{Sec:Method}
In this section, we introduce a fully discrete numerical method for solving \eqref{hodge-varphi-j} and \eqref{PDE-A}-\eqref{Def-A0}, and then present the main theoretical results about the convergence of the numerical solutions. 

Let $\Im_h$, $0<h<1$, be a family of quasi-uniform partition of $\Omega$ into triangles $K_j$, $j=1,\ldots,M$, such that no vertex of any triangle lies on the interior of an edge of another triangle. Here, quasi-uniformity means that the diameter $h_{K_j}$ of a triangle $K_j$ and the radius $\rho_{K_j}$ of the inscribed circle satisfy
$$
c^{-1} h\le \rho_{K_j} \le h_{K_j} \le c\, h,\quad j=1,\dots,M,
$$
where $\displaystyle h=\max_{1\leq j\leq M}\{\mbox{diam}(K_j)\}$ is the mesh size, and $c$ is a constant independent of $h$. 

For any integer $r\geq 1$, we define the Taylor--Hood finite element space $\mathring {\bf S}_h^{r+1}\times S_h^r/\R$ with 
\begin{eqnarray*}
&&
S_h^{r}=\{\chi_h\in H^1(\Omega):
{\chi_h}|_{ K_j}\in P_{r}(K_j),\ \forall K_j\in\Im_h\} , \\
&&S_h^r/\R=\{v\in S_h^r: \mbox{$\int_\Omega v \d x=0$}\} , \\
&&
\mathring S_h^{r+1}=S_h^{r+1}\cap H^1_0(\Omega)  ,\\
&&
\mathring {\bf S}_h^{r+1}=\mathring S_h^{r+1}\times \mathring S_h^{r+1} , 
\end{eqnarray*}where $P_r(K_j)$ is the space of
polynomials of degree $r$ on the triangle $K_j$. 
 
For any given $j=1,\dots,m$, let $\varphi_{j,h}\in S_h^{r+1}$ be the finite element solution of 
\eqref{hodge-varphi-j}, i.e., 
\begin{align}\label{numer-varphi}
(\nabla\varphi_{j,h},\nabla v_h)=0 \quad\forall\,v_h\in \mathring S_h^{r+1}
\end{align}
such that $\varphi_{j,h}=1$ on $\Gamma_j$ and $\varphi_{j,h}=0$ on $\partial\Omega\backslash\Gamma_j$. Let $\beta_{j,h}$, $j=1,\dots,m$, be the constants (independent of space and time) determined by the equations 
\begin{align}\label{betajh}
& \sum_{j=1}^m  \beta_{j,h}  (\nabla\times\varphi_{j,h},\nabla\times\varphi_{i,h}) 
= (\H_0,\nabla\times\varphi_{i,h}) \quad \mbox{for }i=1,\dots,m.
\end{align}

Let $\{t_n=n\tau\}_{n=0}^N$ denote a uniform partition of the time interval $[0,T]$, 
with a step size $\tau=T/N$, and $u^n=u(x,t_n)$. 
A fully discrete numerical scheme for the system \eqref{PDE-A}-\eqref{PDE-u-22} is to find $A_h^n\in \mathring S_h^{r+1}$, $\u_h^n\in \mathring {\bf S}_h^{r+1}$, 
and $p_h^n\in S_h^{r}/\R$ such that 
\begin{align}
&\bigg(\mu\frac{A_h^n-A_h^{n-1}}{\tau} ,a_h\bigg)
+(\sigma^{-1}\nabla A_h^n,\nabla a_h)
- (\mu{\bm u}_h^{n}\times\big(\nabla\times A_h^{n-1}+ \mbox{$\sum_{j=1}^m$}\beta_{j,h} \nabla\times\varphi_{j,h} \big),a_h) \nonumber \\
&
=(\sigma^{-1} J^n ,a_h)   \label{FEM-A} \\[8pt]
&\bigg(\frac{{\bm u}_h^n-\u_h^{n-1}}{\tau} ,{\bm v}_h\bigg)
+\frac12({\bm u}_h^{n-1}\cdot\nabla{\bm u}_h^n,{\bm v}_h)
-\frac12({\bm u}_h^{n-1}\cdot\nabla {\bm v}_h,{\bm u}_h^n)
+(\nu\nabla{\bm u}_h^n,\nabla {\bm v}_h)
-(p_h^n ,\nabla\cdot{\bm v}_h) \nonumber \\
&=({\bm f}^n,{\bm v}_h)  + (\mu\big(\nabla\times A_h^{n-1}+ \mbox{$\sum_{j=1}^m$}\beta_{j,h} \nabla\times\varphi_{j,h} \big)\times \Delta_h A_h^n,{\bm v}_h)  
\label{FEM-u} \\[8pt]
&(\nabla\cdot{\bm u}_h^n,q_h)=0  \label{FEM-p}
\end{align}
hold for all test functions $a_h\in \mathring S_h^{r+1}$, ${\bm v}_h\in \mathring{\bf S}_h^{r+1}$,  $q_h\in S_h^{r}/\R$, and $n=1,2,\dots,N$. 
{\color{blue}
Then the discrete magnetic field can be obtained by 
$${\bf H}_h^n=\nabla\times (A_h^n+\sum_{j=1}^m\beta_{j,h} \varphi_{j,h}) .$$  As a result, ${\bf H}_h^n$ satisfies $\nabla\cdot {\bf H}_h^n=0$ automatically.
}
Here, 
the operator $\Delta_h:\mathring S_h^{r+1}\rightarrow \mathring S_h^{r+1}$ is defined via the duality:
\begin{align}
(\Delta_h A_h^n,a_h)=-(\nabla A_h^n,\nabla a_h)\quad \forall a_h\in \mathring S_h^{r+1}.
\end{align} 
The initial condition $A^0_h\in S_h^{r+1}$ can be determined by 
\begin{align}\label{Eq-A0h}
(\nabla A^0_h,\nabla a_h)=(\H_0,\nabla\times a_h),\quad\forall\, a_h\in \mathring S_h^{r+1} .
\end{align} 
The initial condition for velocity is given by ${\bm u}_h^0=P_h{\bm u}_0$, where $P_h$ denotes the $L^2$ projection from $L^2(\Omega)$ onto $\mathring S_h^{r+1}$, defined by 
\begin{align}\label{L2-projection}
\begin{aligned}
&(P_h\varphi,\chi_h)=(\varphi,\chi_h) && \forall \varphi\in L^2(\Omega),\,\,\, \forall \chi_h\in \mathring S_h^{r+1},  
\end{aligned}
\end{align}
which automatically extends to the vector-valued space $\mathring {\bf S}_h^{r+1}$, i.e., 
\begin{align}\label{L2-projection2}
\begin{aligned}
&(P_h{\bm\varphi},{\bm\chi}_h)=({\bm\varphi},{\bm\chi}_h) && \forall \varphi\in L^2(\Omega),\,\,\,  \forall {\bm\chi}_h\in \mathring {\bf S}_h^{r+1}.  
\end{aligned}
\end{align}

For any sequence $\omega_h^{n}$, $n=1,2,\dots,$ we define the piecewise constant functions $\omega_{h,\tau}^+$ and $\omega_{h,\tau}^-$ by 
\begin{align}
\omega_{h,\tau}^+(t)
:=\omega_h^{n}  
\quad\mbox{and}\quad
\omega_{h,\tau}^-(t)
:=\omega_h^{n-1}  
\end{align}for $t\in(t_{n-1},t_{n}]$ and $n=1,2,\dots,N$. 
Then we have the following result on the convergence of the numerical solutions.\medskip

\begin{theorem}\label{MainTHM}
Under assumption \eqref{Reg-f-J}, the fully discrete finite element scheme \eqref{numer-varphi}-\eqref{FEM-p} has a unique solution. For any sequence $(\tau_n,h_n)\rightarrow (0,0)$, there exists a subsequence $(\tau_{n_k},h_{n_k})\rightarrow (0,0)$ such that the corresponding numerical solutions converge to a weak solution $(A,{\bm u},(\varphi_j)_{j=1}^m,(\beta_j)_{j=1}^m)$ of \eqref{hodge-varphi-j} and \eqref{PDE-A}-\eqref{Def-A0} in the following sense:   
\begin{align*}
&A_{h_{n_k},\tau_{n_k}}^+\,\,\,\mbox{converges to $A$ in}\,\,\,L^2(0,T;H^{1}(\Omega)), \\
&\u_{h_{n_k},\tau_{n_k}}^+\,\,\,\mbox{converges to $\u$ in}\,\,\,L^2(0,T;{\bf L}^2(\Omega)) ,\\
&\varphi_{j,h}\,\,\,\mbox{converges to $\varphi_{j}$ in $H^{1}(\Omega)$},\\
&\beta_{j,h}\,\,\,\mbox{converges to $\beta_{j}$ in $\R$},
\end{align*}
which also imply that 
${\bm H}_{h_{n_k},\tau_{n_k}}^+=
\nabla\times A_{h_{n_k},\tau_{n_k}}^++ \mbox{$\sum_{j=1}^m$}\beta_{j,h} \nabla\times\varphi_{j,h}$ converges to ${\bm H}=\nabla\times A+ \mbox{$\sum_{j=1}^m$}\beta_{j} \nabla\times\varphi_{j}$ in $L^2(0,T;{\bf L}^2(\Omega))$.\medskip
\end{theorem}

\begin{remark}\label{Remark-THM}
{\upshape 
The uniqueness of weak solutions was proved in \cite{Sermange-Temam-1983} for convex and smooth domains, but remains open for nonconvex and nonsmooth domains. Our proof shows that every sequence of numerical solutions contains a subsequence converging to a weak solution of the PDE problem. If the weak solution is unique, then Theorem \ref{MainTHM} implies that the numerical solutions converge to the unique weak solution as $(\tau,h)\rightarrow (0,0)$ (without taking a subsequence). 
{\color{blue}
Theorem \ref{MainTHM} implies the existence of weak solutions for \eqref{PDE-A}-\eqref{Def-A0} with regularity \eqref{RegA}-\eqref{betaj2}. However, the regularity \eqref{RegA} only implies that the right-hand side of \eqref{PDE-u-11} is in $L^2(0,T;{\bf L}^1(\Omega))\hookrightarrow L^2(0,T;{\bf W}^{-1,p}(\Omega))$ for all $1\le p<2$. Whether the right-hand side of \eqref{PDE-u-11} is in $L^2(0,T;{\bf H}^{-1}(\Omega))$ is unknown and requires further investigation. 
}
 
}
\end{remark}

\section{Proof of Theorem \ref{MainTHM}}
In this section, we prove Theorem \ref{MainTHM} by using a compactness argument. 
We first introduce some standard notations of finite element spaces in Section \ref{sec-prel}, and then present energy estimates for the numerical solutions in Section \ref{subsec:Energy}. In Section \ref{sec-compact} we utilize the compactness of the numerical solutions to prove the existence of a subsequence (in every sequence of numerical solutions) that converges to a weak solution of the PDE problem. 

Throughout this paper, we denote by $C$ a generic positive constant which could be different
at different places but would be independent of $n$, $h$, and $\tau$. 
To simplify notation, we  use the abbreviations $W^{k,p}=W^{k,p}(\Omega)$, $L^{p}=L^{p}(\Omega)$ and $H^k=H^{k}(\Omega)$ for $k\ge 0$ and $1\le p\le \infty$.

\subsection{Prelimiaries}\label{sec-prel}

It is known that the Taylor--Hood finite element space $\mathring {\bf S}_h^{r+1}\times S_h^r/\R$ ($r\ge 1$) satisfies the following discrete inf-sup condition for some constant $\gamma >0$: 
\begin{eqnarray}\label{LBB-d}
\sup_{{\bm \chi}_h\in \mathring {\bf S}_h^{r+1}}\frac{(\varphi_h,\nabla\cdot {\bm \chi}_h)}{\|\nabla {\bm \chi}_h\|_{L^2(\Omega)}} \geq
\gamma \|\varphi_h\|_{L^2(\Omega)} \quad \forall \varphi_h\in S_h^r .
\end{eqnarray}

Over the finite element spaces $\mathring S_h^{r+1}$ and $\mathring {\bf S}_h^{r+1}$, we define the $L^2$ projection $P_h$ in \eqref{L2-projection}-\eqref{L2-projection2}. 
Besides, we also define the $L^2$ projection $\widetilde P_h:L^2(\Omega)\rightarrow S_h^r$ (without enforcing boundary conditions), i.e., 
\begin{align*}
\begin{aligned}
&(\widetilde P_h\varphi,\chi_h)=(\varphi,\chi_h) && \forall \varphi\in L^2(\Omega),\,\,\, \forall \chi_h\in S_h ^r . 
\end{aligned}
\end{align*}
We denote by $W^{-1,p}(\Omega)$ the dual space of $W^{-1,p'}_0(\Omega)$. 
Under the assumptions on the triangulation in Section \ref{Sec:Method}, the $L^2$ projections defined above satisfy the following standard estimates for $1\le p\le \infty$:
\begin{align}
&\|P_h\varphi\|_{W^{1,p}(\Omega)}  \le C\|\varphi\|_{W^{1,p}(\Omega)}  &&\forall\,\varphi\in W^{1,p}_0(\Omega) \label{StabPh1}\\
&\|P_h\varphi\|_{L^p(\Omega)}  \le C\|\varphi\|_{L^p(\Omega)}  &&\forall\,\varphi\in L^p(\Omega) \label{StabPh2} \\
&\|P_h\varphi\|_{W^{-1,p}(\Omega)}  \le C\|\varphi\|_{W^{-1,p}(\Omega)}  &&\forall\,\varphi\in W^{-1,p}(\Omega) \label{StabPh3}\\
&\|\widetilde P_h\varphi\|_{L^p(\Omega)}  \le C\|\varphi\|_{L^p(\Omega)}  &&\forall\,\varphi\in L^p(\Omega)  \label{StabPh4}\\
&\lim_{h\rightarrow 0}\|P_h\varphi-\varphi\|_{L^p(\Omega)} \rightarrow 0 &&\forall\,\varphi\in L^p(\Omega) \label{Ph-Lp1}\\
&\lim_{h\rightarrow 0}\|P_h\varphi-\varphi\|_{W^{1,p}(\Omega)} \rightarrow 0 &&\forall\,\varphi\in W^{1,p}_0(\Omega)  \label{Ph-W1p0} \\
&\lim_{h\rightarrow 0} \|\widetilde P_h\varphi-\varphi\|_{L^p(\Omega)} \rightarrow 0 &&\forall\,\varphi\in L^p(\Omega) . \label{Ph-Lp2} 
\end{align}
The $L^p$ and $W^{1,p}$ estimates above were proved in \cite[Section 2]{Crouzeix-Thomee-1987} for meshes more general than quasi-uniformity assumed in this paper. 
Since $W^{-1,p}(\Omega)$ is the dual space of $W^{-1,p'}_0(\Omega)$, the self-adjointness of the $L^2$-projection and \eqref{StabPh1} together imply \eqref{StabPh3}, which is used in the proof for the convergence of numerical solutions; see the last inequality of \eqref{use-W-1p-Ph}. 

For any ${\bm v}\in {\bf H}^1_{0,{\rm div}}(\Omega)$ we define $Q_h{\bm v}$ to be the Fortin projection of ${\bm v}$ onto $S_h^{r+1}$, satisfying 
\begin{align}
\begin{aligned}
(\nabla\cdot({\bm v}-Q_h{\bm v}),q_h)=0 &&\forall\, q_h\in S_h^r
\label{v-vh-Fortin}
\end{aligned}
\end{align}
which has the following property (cf. \cite{GS03} and \cite[Lemma 3.4]{GHL-2018})
\begin{align}
&
\|Q_h{\bm v}\|_{H^1 }
\le C\|{\bm v}\|_{H^1 } &&\forall\, {\bm v}\in {\bf H}^1_{0,{\rm div}}(\Omega)\\
&
\|{\bm v}-Q_h{\bm v}\|_{L^2 }
\le C\|{\bm v}\|_{H^1 }h &&\forall\, {\bm v}\in {\bf H}^1_{0,{\rm div}}(\Omega).
\end{align}

Furthermore, over the finite element spaces $\mathring S_h^{r+1}$, the following inverse inequality holds; see \cite{DFJ-1974,Wheeler1973}.
\begin{align}
&
\|\chi_h\|_{W^{1,q}}
\le
Ch^{-1}\|\chi_h\|_{L^{q}} \\
&
\|\chi_h\|_{W^{m,q}}
\le
Ch^{2/q-2/{\tilde q}}\|\chi_h\|_{W^{m,\tilde q}} 
\end{align}
for all $\chi_h\in \mathring S_h^{r+1}$, and $1\le \tilde q\le q\le \infty$, $m=0,1$.

\subsection{Energy estimate}\label{subsec:Energy}
For any sequence of functions $\omega_h^{n}$, $n=0,1\dots,N$, we let $\omega_{h,\tau}$ denote a piecewise linear function in time defined by
\begin{align}
\omega_{h,\tau}(t)
:=\frac{t_{n}-t}{\tau}\omega_h^{n-1}+\frac{t-t_{n-1}}{\tau}\omega_h^{n}  
\end{align}for $t\in(t_{n-1},t_{n}]$ and $n=1,2,\dots,N$. 
Then we have the following energy estimate for the numerical solutions. 
\begin{proposition}\label{thm-energy}
The numerical scheme \eqref{numer-varphi}-\eqref{FEM-p} admits a unique solution $(A_h^n, {\bm u}_h^n,p_h^n)$, which satisfies the following estimates for $s\in(1,2)$: 
\begin{align}\label{Ahtau-uhtau}
&\|A_{h,\tau}\|_{C([0,T];H^1(\Omega))}
+\|\Delta_h A_{h,\tau}\|_{L^2(0,T;L^2(\Omega))}
+\|\partial_tA_{h,\tau}\|_{L^s(0,T;L^2(\Omega))} 
\nonumber \\
&+\|\u_{h,\tau}\|_{C([0,T];{\bf L}^2(\Omega))} 
+\|\u_{h,\tau}\|_{L^2(0,T;{\bf H}^1(\Omega))} 
+\|\partial_t\u_{h,\tau}\|_{L^s(0,T;{\bf H}^1_{0,{\rm div}}(\Omega)')} 
\le C .
\end{align} 
Furthermore, the finite element solution $(A_{h,\tau}^\pm,{\bm u}_{h,\tau}^\pm)$ satisfies the estimate
\begin{align}\label{Ahtau-uhtau-2}
&\|A_{h,\tau}^\pm\|_{L^\infty(0,T;H^1(\Omega))}
+\|\Delta_h A_{h,\tau}^\pm\|_{L^2(0,T;L^2(\Omega))}
+\|\u_{h,\tau}^\pm\|_{L^\infty(0,T;{\bf L}^2(\Omega))} 
+\|\u_{h,\tau}^\pm\|_{L^2(0,T;{\bf H}^1(\Omega))} 
\le C .
\end{align}
\end{proposition}

To prove Proposition \ref{thm-energy} we need to use the following two lemmas.
\begin{lemma}\label{DeltaA-L4}
There exists $\alpha>0$ (depending only on the domain $\Omega$) such that 
\begin{align}
\|\chi_h\|_{W^{1,4+\alpha}}\le C\|\Delta_h\chi_h\|_{L^2} 
\end{align}
for any $\chi_h\in \mathring S_h^{r+1}$. 
Furthermore, if $h_n \rightarrow 0$ as $n\rightarrow\infty$, 
and $\Delta_h\chi_{h_n}$ is bounded in $L^2(\Omega)$, then 
$\chi_{h_n}$ is compact in $W^{1,4+\alpha}(\Omega)$. 
\end{lemma}

{\it Proof.}$\,\,\,$ 
For any given $\chi_{h_n}\in \mathring S_{h_n}^{r+1}$, let $\chi_n\in H^1_0(\Omega)$ be the weak solution of the PDE problem
\begin{align}\label{PDE-chin}
\Delta \chi_n=\Delta_h\chi_{h_n}
\end{align}
with homogeneous Dirichlet boundary condition. Thus $\chi_n$ satisfies the weak formulation
\begin{align*}
(\nabla \chi_n,\nabla\chi)
=-(\Delta_h\chi_{h_n},\chi)
=(\nabla \chi_{h_n},\nabla\chi) 
\end{align*}for any $\chi\in \mathring S_{h_n}^{r+1}$, 
which implies $\chi_{h_n}$ is the Ritz projection of $\chi_n$. 

On the one hand, as the solution of the PDE problem \eqref{PDE-chin}, $\chi_n$ satisfies the standard PDE estimate
\begin{align}\label{Am-H32}
\|\chi_n\|_{H^{\frac32+\delta}} 
\leq C\|\Delta_h\chi_{h_n}\|_{L^2} 
\end{align}
for some constant $\delta\in(0,\frac{1}{2}]$ depending on the maximal interior angle of the domain $\Omega$. 
The estimate \eqref{Am-H32} is a consequence of \cite[Corollary 3.9, with fractional $k$]{Dauge}; also see \cite[p. 30]{Dauge08} and \cite[(23.3)]{Dauge88}. 
Since $H^{\frac32+\delta}(\Omega)$ is compactly embedded into $W^{1,4+\alpha}(\Omega)$ for $\alpha\in(0,\frac{8\delta}{1-2\delta})$ in a polygon (cf. \cite[Theorem 7.34]{Adams-Fournier-2003}, with $p=q=2$ and $s=\frac12+\delta$), it follows that 
\begin{align}\label{eq-chin}
\|\chi_n\|_{W^{1,4+\alpha}}
\le
C\|\chi_n\|_{H^{\frac32+\delta}} 
\leq C\|\Delta_h\chi_{h_n}\|_{L^2} 
\end{align}
and the set of functions $\{\chi_n:n=1,2,\dots\}$ is compact in $W^{1,4+\alpha}(\Omega)$. 

On the other hand, as the Ritz projection of $\chi_n$, the finite element function $\chi_{h_n}$ satisfies the standard error estimate
\begin{align*}
\|\chi_{h_n}-\chi_n\|_{H^1}
\le C\|\chi_n\|_{H^{\frac32+\delta}} h^{\frac12+\delta} 
\le C\|\Delta_h\chi_{h_n}\|_{L^2} h_n^{\frac12+\delta} \rightarrow 0
\end{align*}as $n\rightarrow \infty$. 
By using the triangle inequality and the inverse inequality of the finite element space, we have 
\begin{align*}
&\|\chi_{h_n}-\chi_n\|_{W^{1,4+\alpha}} \\ 
&\le \|\chi_{h_n}-P_h\chi_n\|_{W^{1,4+\alpha}} 
+\|P_h\chi_n-\chi_n\|_{W^{1,4+\alpha}} \\
&\le Ch_n^{-\frac{2+\alpha}{4+\alpha}}\|\chi_{h_n}-P_h\chi_n\|_{H^1} 
+\|P_h\chi_n-\chi_n\|_{W^{1,4+\alpha}}\\
&\le Ch_n^{-\frac{2+\alpha}{4+\alpha}}\|\chi_{h_n}-\chi_n\|_{H^1}
+Ch_n^{-\frac{2+\alpha}{4+\alpha}}\|\chi_n-P_h\chi_n\|_{H^1} 
+\|P_h\chi_n-\chi_n\|_{W^{1,4+\alpha}}\\
&\le C\|\chi_n\|_{H^{\frac32+\delta}} h_n^{\frac12+\delta-\frac{2+\alpha}{4+\alpha}}
+C\|\chi_n\|_{H^{\frac32+\delta}} h_n^{\frac12+\delta-\frac{2+\alpha}{4+\alpha}} 
+C\|\chi_n\|_{H^{\frac32+\delta}} h_n^{\frac12+\delta-\frac{2+\alpha}{4+\alpha}}   \\
&\le C\|\Delta_h\chi_{h_n}\|_{L^2} h_n^{\frac12+\delta-\frac{2+\alpha}{4+\alpha}}  \\
&\rightarrow 0 \quad\mbox{as $n\rightarrow \infty$, }
\end{align*}
where we have used the fact $\frac12+\delta-\frac{2+\alpha}{4+\alpha}>0$ for $\alpha\in(0,\frac{8\delta}{1-2\delta})$ in the last inequality. 
Thus, by \eqref{eq-chin}, we obtain
\begin{align*}
\|\chi_{h_n}\|_{W^{1,4+\alpha}}
\le
\|\chi_n\|_{W^{1,4+\alpha}}
+C\|\Delta_h\chi_{h_n}\|_{L^2} 
\le
C\|\Delta_h\chi_{h_n}\|_{L^2}.
\end{align*}
Since $\{\chi_n:n=1,2,\dots\}$ is compact in $W^{1,4+\alpha}(\Omega)$ and with the above result, the set of functions $\{\chi_{h_n}:\, n=1,2,\dots\}$ 
is also compact in $W^{1,4+\alpha}(\Omega)$. This completes the proof of Lemma \ref{DeltaA-L4}. 
\hfill\endproof

\begin{lemma}\label{Lemma-varphij}
There exists $\alpha>0$ (depending only on the domain $\Omega$) such that the function $\varphi_{j,h}$ deteremind by \eqref{numer-varphi} converges to $\varphi_j$ in the following sense: 
\begin{align}\label{Est-varphi-j}
&\lim_{h\rightarrow 0}
\|\varphi_{j,h}-\varphi_{j}\|_{W^{1,4+\alpha}(\Omega)} =0 .
\end{align}
Furthermore,
\begin{align}\label{Est-beta-j}
&\lim_{h\rightarrow 0}
|\beta_{j,h}-\beta_{j}| =0 .
\end{align}
\end{lemma}
{\it Proof.}$\,\,$ 
Let $\chi$ be a smooth cut-off function such that $\chi=1$ in a neighborhood of $\Omega_j$ and $\chi=0$ on $\cup_{i\neq j}\partial\Omega_i$. Then $\varphi_j-\chi$ is the solution of 
\begin{align}\label{hodge-varphi-j-2}
\left\{\begin{aligned}
&\Delta (\varphi_j-\chi)=-\Delta \chi &&\mbox{in}\,\,\,\Omega \\
&\varphi_j-\chi=0 &&\mbox{on}\,\,\,\partial\Omega 
\end{aligned}\right.
\end{align}
which implies (similar as \eqref{Am-H32}) 
\begin{align*}
\|\varphi_j-\chi\|_{H^{\frac32+\delta}} 
\leq C\|\Delta \chi\|_{L^2} \le C \quad\mbox{for some}\,\,\,\delta>0\,\,\,\mbox{depending only on the domain $\Omega$}. 
\end{align*}
Therefore, $\|\varphi_j\|_{H^{\frac32+\delta}} \le C$. 
In view of the embedding $H^{\frac32+\delta}(\Omega)\hookrightarrow W^{1,4+\alpha}(\Omega)$ for $\alpha\in(0,\frac{8\delta}{1-2\delta})$ {\color{blue}(we denote by ``$\hookrightarrow$" continuously embedding)}, we have $\nabla\times\varphi_{j}\in L^{4+\alpha}(\Omega)$. 
Since $\varphi_{j,h}$ is the finite element solution of $\varphi_{j}$, it follows that 
\begin{align}
\|P_h\varphi_j-\varphi_{j,h}\|_{H^1} 
\leq Ch^{\frac12+\delta} ,
\end{align}
where $P_h$ is the $L^2$-projection operator, satisfying \eqref{StabPh1}. By using the inverse inequality we obtain 
\begin{align}
\|P_h\varphi_j-\varphi_{j,h}\|_{W^{1,4+\alpha}} 
\le Ch^{-\frac{2+\alpha}{4+\alpha}} \|P_h\varphi_j-\varphi_{j,h}\|_{H^1} 
\le Ch^{\frac{1}{2}+\delta-\frac{2+\alpha}{4+\alpha}} . 
\end{align}
For $\alpha\in(0,\frac{8\delta}{1-2\delta})$ we have $\frac{1}{2}+\delta-\frac{2+\alpha}{4+\alpha}>0$. 
Since 
$H^{\frac32+\delta}(\Omega) \hookrightarrow W^{\frac{3}{2}+\delta-\frac{2+\alpha}{4+\alpha},4+\alpha}(\Omega)$, it follows that 
\begin{align}
\|\varphi_j-P_h\varphi_j\|_{W^{1,4+\alpha}}
\le
C\|\varphi_j\|_{W^{\frac{3}{2}+\delta-\frac{2+\alpha}{4+\alpha},4+\alpha}} h^{\frac{1}{2}+\delta-\frac{2+\alpha}{4+\alpha}}
\le
C\|\varphi_j\|_{H^{\frac{3}{2}+\delta}}h^{\frac{1}{2}+\delta-\frac{2+\alpha}{4+\alpha}}
\le
Ch^{\frac{1}{2}+\delta-\frac{2+\alpha}{4+\alpha}}.
\end{align} 
Threfore we have 
\begin{align}
\|\varphi_j-\varphi_{j,h}\|_{W^{1,4+\alpha}} 
\le \|\varphi_j-P_h\varphi_{j}\|_{W^{1,4+\alpha}} 
+\|P_h\varphi_j-\varphi_{j,h}\|_{W^{1,4+\alpha}}  
\le  Ch^{\frac{1}{2}+\delta-\frac{2+\alpha}{4+\alpha}} 
\end{align}
which implies \eqref{Est-varphi-j}.

Let $M=(M_{ij})$ and $M_h=(M_{ij,h})$ be two $m\times m$ matrices, with $M_{ij}=(\nabla\times\varphi_j,\nabla\times\varphi_i)$ and $M_{ij,h}=(\nabla\times\varphi_{j,h},\nabla\times\varphi_{i,h})$. 
Since $M=(M_{ij})$ is positive definite and \eqref{Est-varphi-j} implies $\lim_{h\rightarrow 0}M_{h}=M$, it follows that 
$\lim_{h\rightarrow 0}M_{h}^{-1}=M^{-1}$. 
Then 
\eqref{betajh} and \eqref{betaj} imply 
\begin{align*}
\lim_{h\rightarrow 0} \beta_{i,h}
&=\lim_{h\rightarrow 0} 
\sum_{j=1}^m (M_h^{-1})_{ij}(\H_0,\nabla\times\varphi_{j,h}) \\ 
&= \sum_{j=1}^m (M^{-1})_{ij} \lim_{h\rightarrow 0}  (\H_0,\nabla\times\varphi_{j,h}) \\
&= \sum_{j=1}^m (M^{-1})_{ij}  (\H_0,\nabla\times\varphi_j) \\
&=\beta_{i} . 
\end{align*}
\hfill\endproof 

With above results, we start to prove Proposition \ref{thm-energy}. 

{\it Proof of Proposition \ref{thm-energy}.}$\,\,\,$ 
First, we start with estimating $\|\nabla A_h^n\|_{L^2}$ and $\|{\bm u}_h^n\|_{L^2}$. 
To this end, we substitute $a_h=-\Delta_hA_h^n$, ${\bm v}_h=\u_h^n$ and $q_h=p_h^n$ into \eqref{FEM-A}-\eqref{FEM-u}. Then we obtain 
\begin{align}
&\bigg(\mu\frac{\nabla A_h^n-\nabla A_h^{n-1}}{\tau} ,\nabla A_h^n\bigg)
+\sigma^{-1}\|\Delta_h A_h^n\|_{L^2}^2 
+ (\mu{\bm u}_h^{n}\times\big(\nabla\times A_h^{n-1}+ \mbox{$\sum_{j=1}^m$}\beta_{j,h} \nabla\times\varphi_{j,h} \big),\Delta_h A_h^n) \nonumber \\
&
=(\sigma^{-1} J^n ,\Delta_h A_h^n)  \label{Test1A} \\
&\bigg(\frac{\u_h^n-\u_h^{n-1}}{\tau} ,\u_h\bigg)
+\nu\|\nabla\u_h^n\|_{L^2}^2 
=({\bm f}^n,\u_h^n)  + (\mu\big(\nabla\times A_h^{n-1}+ \mbox{$\sum_{j=1}^m$}\beta_{j,h} \nabla\times\varphi_{j,h} \big)\times \Delta_h A_h^n,\u_h^n) \nonumber \\
&\hspace{139pt}  =({\bm f}^n,\u_h^n)  + (\mu\u_h^n\times \big(\nabla\times A_h^{n-1}+ \mbox{$\sum_{j=1}^m$}\beta_{j,h} \nabla\times\varphi_{j,h} \big), \Delta_h A_h^n) . \label{Test1u} 
\end{align}
{\color{blue}
Since 
\begin{align*}
&\bigg(\mu\frac{\nabla A_h^n-\nabla A_h^{n-1}}{\tau} ,\nabla A_h^n\bigg)
= \frac{\mu\|\nabla A_h^n\|_{L^2}^2-\mu\|\nabla A_h^{n-1}\|_{L^2}^2}{2\tau} 
+  \frac{\mu\|\nabla A_h^n-\nabla A_h^{n-1}\|_{L^2}^2}{2\tau} ,\\
&\bigg(\frac{\u_h^n-\u_h^{n-1}}{\tau} ,\u_h\bigg)
= \frac{\|\u_h^n\|_{L^2}^2-\|\u_h^{n-1}\|_{L^2}^2}{2\tau}
+  \frac{\|\u_h^n-\u_h^{n-1}\|_{L^2}^2}{2\tau} ,
\end{align*}
summing up \eqref{Test1A}-\eqref{Test1u} and dropping the square terms in the above two equations yield} 
\begin{align}
&\frac{\mu\|\nabla A_h^n\|_{L^2}^2-\mu\|\nabla A_h^{n-1}\|_{L^2}^2}{2\tau}
+\frac{\|\u_h^n\|_{L^2}^2-\|\u_h^{n-1}\|_{L^2}^2}{2\tau}
+\sigma^{-1}\|\Delta_h A_h^n\|_{L^2}^2 
+\nu\|\nabla\u_h^n\|_{L^2}^2 \nonumber \\ 
&\le (\sigma^{-1} J^n ,\Delta_h A_h^n) 
     +({\bm f}^n,\u_h) \nonumber \\
&\le \frac{\sigma^{-1}}{2} \|J^n\|_{L^2}^2+\frac{\sigma^{-1}}{2}\|\Delta_h A_h^n\|_{L^2}^2
     +\frac{1}{4\epsilon} \|{\bm f}^n\|_{L^2}^2+\epsilon\|\u_h^n\|_{L^2}^2 \nonumber \\
     &\le \frac{\sigma^{-1}}{2} \|J^n\|_{L^2}^2+\frac{\sigma^{-1}}{2}\|\Delta_h A_h^n\|_{L^2}^2
     +\frac{1}{4\epsilon} \|{\bm f}^n\|_{L^2}^2+C\epsilon\|\nabla\u_h^n\|_{L^2}^2 .
\end{align}
By choosing a sufficiently small $\epsilon$, the inequality above implies 
\begin{align}\label{im-est}
&\frac{\mu\|\nabla A_h^n\|_{L^2}^2-\mu\|\nabla A_h^{n-1}\|_{L^2}^2}{2\tau}
+\frac{\|\u_h^n\|_{L^2}^2-\|\u_h^{n-1}\|_{L^2}^2}{2\tau}
+\frac{\sigma^{-1}}{2}\|\Delta_h A_h^n\|_{L^2}^2 
+\frac{\nu}{2}\|\nabla\u_h^n\|_{L^2}^2 \nonumber \\
&\le C(\|J^n\|_{L^2}^2 +\|{\bm f}^n\|_{L^2}^2) .
\end{align}
By summing up the inequality above for $n=1,\dots,m$, with $1\le m\le N$, we obtain  
\begin{align}
&\max_{1\le n\le N} (\|\nabla A_h^n\|_{L^2}^2+\|\u_h^n\|_{L^2}^2)
+\tau \sum_{n=1}^N(\|\Delta_h A_h^n\|_{L^2}^2 +\|\nabla\u_h^n\|_{L^2}^2) \nonumber \\
&\le C\tau \sum_{n=1}^N(\|J^n\|_{L^2}^2 +\|{\bm f}^n\|_{L^2}^2) 
+C(\|\nabla A_h^0\|_{L^2}^2+\|\u_h^0\|_{L^2}^2) \nonumber  \\
&\le C ,
\label{FEM-E}  
\end{align}
where we have used the boundedness of $\|\nabla A_h^0\|_{L^2}$ and $\|\u_h^0\|_{L^2}$, which are direct consequences of the definitions of $\u_h^0=P_h\u_0$ and $A_h^0$ in \eqref{Eq-A0h}. 
The estimate \eqref{im-est} also implies that, by setting $A_h^{n-1}=J^n=0$ and $\u_h^{n-1}={\bm f}^n={\bm 0}$, the homogeneous linear system given by \eqref{numer-varphi}-\eqref{FEM-p} has only zero solution $A_h^n=0$ and $\u_h^n=0$. 
{\color{blue}
In this case, \eqref{FEM-u} reduces to 
\begin{align*}
(p_h^n ,\nabla\cdot{\bm v}_h) = 0 .
\end{align*}
By using the inf-sup condition \eqref{LBB-d}, we further derive $p_h^n=0$.} 
This implies the existence and uniqueness of solutions of the linear system \eqref{numer-varphi}-\eqref{FEM-p}. 

Second, we estimate $\big\|\frac{A_h^n-A_h^{n-1}}{\tau}\big\|_{L^2}$. 
To this end, we use \eqref{FEM-A} and note that for any $q,\bar q>2$ satisfying $\frac{1}{\bar q}+\frac{1}{q}=\frac{1}{2}$ 
\begin{align*}
&\bigg|\bigg(\mu\frac{A_h^n-A_h^{n-1}}{\tau} ,a_h\bigg)\bigg| \\ 
&= \bigg|(\sigma^{-1}\Delta_h A_h^n,a_h)
+(\mu{\bm u}_h^{n}\times\big(\nabla\times A_h^{n-1}+ \mbox{$\sum_{j=1}^m$}\beta_{j,h} \nabla\times\varphi_{j,h} \big),a_h)
+(\sigma^{-1} J^n ,a_h) \bigg| \\
&\le (\sigma^{-1}\|\Delta_h A_h^n\|_{L^2}
+\mu\|{\bm u}_h^{n}\times\big(\nabla\times A_h^{n-1}+ \mbox{$\sum_{j=1}^m$}\beta_{j,h} \nabla\times\varphi_{j,h} \big)\|_{L^2}
+\sigma^{-1}\|J^n\|_{L^2})\|a_h\|_{L^2} \\
&\le C(\|\Delta_h A_h^n\|_{L^2}
+\|{\bm u}_h^{n}\|_{L^{\bar q}}\|\nabla\times A_h^{n-1}+ \mbox{$\sum_{j=1}^m$}\beta_{j,h} \nabla\times\varphi_{j,h} \|_{L^q}
+\|J^n\|_{L^2})\|a_h\|_{L^2} \\
&\le C(\|\Delta_h A_h^n\|_{L^2} 
+\|\nabla {\bm u}_h^{n}\|_{L^2}\|\nabla\times A_h^{n-1}+ \mbox{$\sum_{j=1}^m$}\beta_{j,h} \nabla\times\varphi_{j,h} \|_{L^q} 
+\|J^n\|_{L^2})\|a_h\|_{L^2} .
\end{align*}
The inequality above implies 
(via the duality argument)
\begin{align}\label{jjaaa}
\bigg\|\frac{A_h^n-A_h^{n-1}}{\tau}\bigg\|_{L^2}
\le C(\|\Delta_h A_h^n\|_{L^2}
+\|\nabla{\bm u}_h^{n}\|_{L^2}\|\nabla\times A_h^{n-1}+ \mbox{$\sum_{j=1}^m$}\beta_{j,h} \nabla\times\varphi_{j,h} \|_{L^q}
+\|J^n\|_{L^2}) 
\end{align}
which holds for all $q\in(2,4)$ with $\bar q=\frac{2q}{q-2}$. 
By using Lemma \ref{DeltaA-L4} and \eqref{FEM-E}, we have 
\begin{align*}
\|\nabla\times A_h^{n-1}\|_{L^q}\le C\|\nabla A_h^{n-1}\|_{L^2}^{\frac{4}{q}-1}\|\nabla A_h^{n-1}\|_{L^4}^{2-\frac{4}{q}}\le C\|\Delta_h A_h^{n-1}\|_{L^2}^{2-\frac{4}{q}} .
\end{align*}
Therefore, 
\begin{align*}
&\bigg(\tau\sum_{n=1}^N \|\nabla {\bm u}_h^{n}\|_{L^2}^s \|\nabla\times A_h^{n-1}+ \mbox{$\sum_{j=1}^m$}\beta_{j,h} \nabla\times\varphi_{j,h} \|_{L^q}^s\bigg)^{\frac1s} \nonumber \\
&\le 
C\bigg(\tau\sum_{n=1}^N \|\nabla {\bm u}_h^{n}\|_{L^2}^2\bigg)^{\frac12}
\bigg(\tau\sum_{n=1}^N \|\nabla\times A_h^{n-1}+ \mbox{$\sum_{j=1}^m$}\beta_{j,h} \nabla\times\varphi_{j,h} \|_{L^q}^{\frac{2s}{2-s}}\bigg)^{\frac{2-s}{2s}}\\
&\le 
C\bigg(\tau\sum_{n=1}^N \|\nabla {\bm u}_h^{n}\|_{L^2}^2\bigg)^{\frac12}
\bigg(\tau\sum_{n=1}^N\|\Delta_h A_h^{n-1}\|_{L^2}^{(2-\frac{4}{q})\frac{2s}{2-s}} + C \bigg)^{\frac{2-s}{2s}} , 
\end{align*}
where we have used Lemma \ref{Lemma-varphij} in the last inequality. 
For any given $1<s<2$ one can choose $q=\frac{4s}{3s-2}\in(2,4)$ so that $(2-\frac{4}{q})\frac{2s}{2-s}=2$ and 
\begin{align*}
&\bigg(\tau\sum_{n=1}^N \|\nabla {\bm u}_h^{n}\|_{L^2}^s \|\nabla\times A_h^{n-1}+ \mbox{$\sum_{j=1}^m$}\beta_{j,h} \nabla\times\varphi_{j,h} \|_{L^q}^s\bigg)^{\frac1s} \nonumber \\
&\le 
C\bigg(\tau\sum_{n=1}^N \|\nabla {\bm u}_h^{n}\|_{L^2}^2\bigg)^{\frac12}
\bigg(\tau\sum_{n=1}^N\|\Delta_h A_h^{n-1}\|_{L^2}^2+C\bigg)^{\frac{2-s}{2s}} 
\nonumber \\
&\le C ,
\end{align*}
where we have used \eqref{FEM-E} again. 
Substituting the estimate above into \eqref{jjaaa} yields 
\begin{align}\label{DtauA-L2Lp}
\bigg(\tau\sum_{n=1}^N \bigg\|\frac{A_h^n-A_h^{n-1}}{\tau}\bigg\|_{L^2}^s\bigg)^{\frac1s} 
\le C ,\quad\forall\, s\in(1,2). 
\end{align}

Third, we present the estimate for $\big\|\frac{{\bm u}_h^n-{\bm u}_h^{n-1}}{\tau}\big\|_{{\bf H}^1_{0,{\rm div}}(\Omega)'}$.
From \eqref{FEM-u} we can derive that 
\begin{align}\label{Dtau-uhn-g}
\bigg(\frac{{\bm u}_h^n-\u_h^{n-1}}{\tau} ,{\bm v}_h\bigg)
&=
(P_h{\bm g}^n,{\bm v}_h)+(p_h^n ,\nabla\cdot{\bm v}_h) , 
\end{align}
where
\begin{align*}
{\bm g}^n
=&-\frac12 {\bm u}_h^{n-1}\cdot\nabla{\bm u}_h^n
-\frac12 (\nabla\cdot{\bm u}_h^{n-1})\cdot{\bm u}_h^n
-\frac12 {\bm u}_h^{n-1}\cdot\nabla {\bm u}_h^n 
-\nu\Delta_h{\bm u}_h^n 
 \nonumber \\
&+ {\bm f}^n  + \mu(\nabla\times A_h^{n-1}+ \mbox{$\sum_{j=1}^m$}\beta_{j,h} \nabla\times\varphi_{j,h})\times \Delta_h A_h^n .
\end{align*}
Substituting ${\bm v}_h=\frac{{\bm u}_h^n-\u_h^{n-1}}{\tau}$ into \eqref{Dtau-uhn-g} yields 
\begin{align*}
\bigg\|\frac{{\bm u}_h^n-\u_h^{n-1}}{\tau} \bigg\|_{L^2}^2 
&\le C\|P_h{\bm g}^n\|_{L^2} \bigg\|\frac{{\bm u}_h^n-\u_h^{n-1}}{\tau} \bigg\|_{L^2}+\bigg(p_h^n ,\frac{\nabla\cdot {\bm u}_h^n-\nabla\cdot \u_h^{n-1}}{\tau} \bigg) \\
&=C\|P_h{\bm g}^n\|_{L^2} \bigg\|\frac{{\bm u}_h^n-\u_h^{n-1}}{\tau} \bigg\|_{L^2}, 
\end{align*}
where we have used \eqref{FEM-p}. 
Then, 
\begin{align}\label{Dtau-L2}
\bigg\|\frac{{\bm u}_h^n-\u_h^{n-1}}{\tau} \bigg\|_{L^2} 
&\le C\|P_h{\bm g}^n\| _{L^2} . 
\end{align}
From \eqref{v-vh-Fortin}, we have $(\nabla\cdot{\bm v},p_h^n)=(\nabla\cdot Q_h{\bm v},p_h^n)=0$ for any ${\bm v}\in {\bf H}^1_{0,{\rm div}}(\Omega)$ and thus 
\begin{align}\label{Dtau-uh-v}
&\bigg(\frac{{\bm u}_h^n-\u_h^{n-1}}{\tau} ,{\bm v}\bigg) \nonumber \\ 
&=\bigg(\frac{{\bm u}_h^n-\u_h^{n-1}}{\tau} ,Q_h{\bm v}\bigg)
+\bigg(\frac{{\bm u}_h^n-\u_h^{n-1}}{\tau} ,{\bm v}-Q_h{\bm v}\bigg) \nonumber \\ 
&=(P_h{\bm g}^n,Q_h{\bm v}) 
     +\bigg(\frac{{\bm u}_h^n-\u_h^{n-1}}{\tau} ,{\bm v}-Q_h{\bm v}\bigg) \nonumber \\ 
&\le C\|{\bm g}^n\|_{H^{-1}(\Omega)}\|{\bm v}\|_{H^1}
+C\bigg\|\frac{{\bm u}_h^n-\u_h^{n-1}}{\tau}\bigg\|_{L^2}\|{\bm v}-Q_h{\bm v}\|_{L^2}   ,
\quad\forall\, \v\in {\bf H}^1_{0,{\rm div}}(\Omega). 
\end{align}
By using H\"older's inequality and \eqref{v-vh-Fortin}, we have 
\begin{align}\label{use-W-1p-Ph}
\bigg\|\frac{{\bm u}_h^n-\u_h^{n-1}}{\tau}\bigg\|_{L^2}\|{\bm v}-Q_h{\bm v}\|_{L^2}
&\le Ch\bigg\|\frac{{\bm u}_h^n-\u_h^{n-1}}{\tau}\bigg\|_{L^2 } \|{\bm v}\|_{H^1} \nonumber \\ 
&\le Ch \|P_h{\bm g}^n\|_{L^2 } \|{\bm v}\|_{H^1} \nonumber \\
&\le C \|P_h{\bm g}^n\|_{H^{-1} } \|{\bm v}\|_{H^1 } \nonumber \\
&\le C \|{\bm g}^n\|_{H^{-1} } \|{\bm v}\|_{H^1 } ,
\end{align}
where we have used the inverse inequality in the second to last inequality, and \eqref{StabPh3} in the last inequality. The two estimates above imply
\begin{align}\label{Dtau-uh-g}
\bigg\|\frac{{\bm u}_h^n-\u_h^{n-1}}{\tau}\bigg\|_{{\bf H}^1_{0,{\rm div}}(\Omega)'}
\le C\|{\bm g}^n\|_{H^{-1} } .
\end{align}
From \eqref{Dtau-uhn-g} we obtain for any $q>1$ 
\begin{align}\label{H-1-g-0}
&\|{\bm g}^n\|_{H^{-1} } \nonumber \\ 
&\le
C\|{\bm u}_h^{n-1}\cdot\nabla{\bm u}_h^n\|_{H^{-1}}
+C\hspace{-10pt}\sup_{{\bm v}_h\in \mathring {\bf S}_h^{r+1},{\bm v}_h\neq {\bm 0}}\hspace{-10pt}\frac{((\nabla\cdot{\bm u}_h^{n-1})\cdot{\bm u}_h^n+{\bm u}_h^{n-1}\cdot\nabla{\bm u}_h^n,{\bm v}_h)}{\|{\bm v}_h\|_{H_0^1}}  
\nonumber\\
&\quad+C\hspace{-10pt}\sup_{{\bm v}_h\in \mathring {\bf S}_h^{r+1},{\bm v}_h\neq {\bm 0}}\hspace{-10pt}\frac{(\Delta_h{\bm u}_h^n,{\bm v}_h)}{\|{\bm v}_h\|_{H_0^1}} +C\|{\bm f}^n\|_{H^{-1}} +C\|(\nabla\times A_h^{n-1}+ \mbox{$\sum_{j=1}^m$}\beta_{j,h} \nabla\times\varphi_{j,h})\times\Delta_h A_h^n\|_{H^{-1}} 
\nonumber\\
&=
C\|{\bm u}_h^{n-1}\cdot\nabla{\bm u}_h^n\|_{H^{-1}}
+C\hspace{-10pt}\sup_{{\bm v}_h\in \mathring {\bf S}_h^{r+1},{\bm v}_h\neq {\bm 0}}\hspace{-10pt}\frac{({\bm u}_h^n,{\bm u}_h^{n-1}\cdot\nabla{\bm v}_h)}{\|{\bm v}_h\|_{H_0^1}}  
\nonumber\\
&\quad+C\hspace{-10pt}\sup_{{\bm v}_h\in \mathring {\bf S}_h^{r+1},{\bm v}_h\neq {\bm 0}}\hspace{-10pt}\frac{(\nabla{\bm u}_h^n,\nabla{\bm v}_h)}{\|{\bm v}_h\|_{H_0^1}} +C\|{\bm f}^n\|_{H^{-1}} +C\|(\nabla\times A_h^{n-1}+ \mbox{$\sum_{j=1}^m$}\beta_{j,h} \nabla\times\varphi_{j,h})\times\Delta_h A_h^n\|_{H^{-1}} 
\nonumber\\
&\le
C\|{\bm u}_h^{n-1}\cdot\nabla{\bm u}_h^n\|_{L^q}
+C\|{\bm u}_h^{n-1}\|_{L^4}\|{\bm u}_h^n\|_{L^4} \nonumber \\
&\quad 
+C\|\nabla{\bm u}_h^n\|_{L^2}
+C\|{\bm f}^n\|_{H^{-1}} +C\|(\nabla\times A_h^{n-1}+ \mbox{$\sum_{j=1}^m$}\beta_{j,h} \nabla\times\varphi_{j,h})\times\Delta_h A_h^n\|_{L^q} \nonumber\\ 
&\le C\|{\bm u}_h^{n-1}\|_{L^{\frac{2q}{2-q}}}\|\nabla{\bm u}_h^n\|_{L^2}
+C\|{\bm u}_h^{n-1}\|_{L^4 }\|{\bm u}_h^n\|_{L^4 }   \nonumber \\
&\quad +C\|{\bm u}_h^n\|_{H^1 }
+C\|{\bm f}^n\|_{H^{-1} }+C(\|\nabla A_h^{n-1}\|_{L^{\frac{2q}{2-q}}}
+\|\nabla \varphi_{j,h}\|_{L^{\frac{2q}{2-q}}})\|\Delta_hA_h^n\|_{L^2} ,\end{align}
where we have used Sobolev embeddding $L^q(\Omega)\hookrightarrow H^{-1}(\Omega)$ ($\forall\,q>1$) in the second to last inequality, and H\"older's inequality with 
$\frac{1}{q}=\frac{1}{2}+\frac{1}{2q/(2-q)}$ in the last inequality. 
In particular, for $1<q<\frac43$ we have $\frac{2q}{2-q}<4$ and so 
\begin{align*}
&\|{\bm u}_h^{n-1}\|_{L^{\frac{2q}{2-q}}}
\le C\|{\bm u}_h^{n-1}\|_{L^2}^{\frac{2}{q}-1}\|\nabla{\bm u}_h^{n-1}\|_{L^2}^{2-\frac{2}{q}} \\
&\|\nabla A_h^{n-1}\|_{L^{\frac{2q}{2-q}}}
\le \|\nabla A_h^{n-1}\|_{L^2 }^{\frac4q-3}\|\nabla A_h^{n-1}\|_{L^4}^{4-\frac4q} \\
&\|\nabla \varphi_{j,h}^{n-1}\|_{L^{\frac{2q}{2-q}}}\le C. 
\end{align*}
Substituting the three inequalities above into \eqref{H-1-g-0} yields 
\begin{align}\label{H-1-g}
&\|{\bm g}^n\|_{H^{-1} } \nonumber \\ 
&\le C\|{\bm u}_h^{n-1}\|_{L^2}^{\frac{2}{q}-1}
\|\nabla{\bm u}_h^{n-1}\|_{L^2}^{2-\frac{2}{q}}\|\nabla{\bm u}_h^n\|_{L^2 } 
+C\|{\bm u}_h^{n-1}\|_{L^2}^{\frac12}\|\nabla{\bm u}_h^{n-1}\|_{L^2 }^{\frac12}
\|{\bm u}_h^{n}\|_{L^2 }^{\frac12}\|\nabla{\bm u}_h^{n}\|_{L^2 }^{\frac12} \nonumber \\
&\quad +C\|{\bm u}_h^n\|_{H^1 }
+C\|{\bm f}^n\|_{H^{-1} }
+(C\|\nabla A_h^{n-1}\|_{L^2 }^{\frac4q-3}\|\nabla A_h^{n-1}\|_{L^4}^{4-\frac4q}+C)\|\Delta_hA_h^n\|_{L^2 } \nonumber \\
&\le C\|\nabla{\bm u}_h^{n-1}\|_{L^2 }^{2-\frac2q}\|\nabla{\bm u}_h^n\|_{L^2 } 
+C\|\nabla{\bm u}_h^{n-1}\|_{L^2 }^{\frac12}
\|\nabla{\bm u}_h^{n}\|_{L^2 }^{\frac12} \nonumber \\
&\quad +C\|{\bm u}_h^n\|_{H^1 }
+C\|{\bm f}^n\|_{H^{-1} }
+(C\|\Delta_hA_h^{n-1}\|_{L^2}^{4-\frac4q} +C)\|\Delta_hA_h^n\|_{L^2} \nonumber\\
&\le C\|\nabla{\bm u}_h^{n-1}\|_{L^2 }^{2-\frac2q}\|\nabla{\bm u}_h^n\|_{L^2 } 
+C\|\nabla{\bm u}_h^{n-1}\|_{L^2 }^{\frac12}
\|\nabla{\bm u}_h^{n}\|_{L^2 }^{\frac12} \nonumber \\
&\quad +C\|{\bm u}_h^n\|_{H^1 }
+C\|{\bm f}^n\|_{H^{-1} }
+C\|\Delta_hA_h^{n-1}\|_{L^2}^{5-\frac4q} 
+C\|\Delta_hA_h^{n}\|_{L^2}^{5-\frac4q} +C 
\end{align}
for all $q\in(1,4/3)$, 
where we have used \eqref{FEM-E}
in the second to last inequality and Young's inequality at the last step. 
For any $s\in(1,2)$ we can choose $q$ to be sufficiently close to $1$ so that $s(3-\frac2q)\le 2$ and $s(5-\frac4q)\le 2$, and therefore 
\begin{align}\label{Dtau-g-1}
\tau\sum_{n=1}^N\|{\bm g}^n\|_{H^{-1} }^s 
&\le 
C\tau\sum_{n=1}^N
\Big(
\|\nabla{\bm u}_h^{n-1}\|_{L^2}^{s(2-\frac2q)}\|\nabla{\bm u}_h^n\|_{L^2}^s 
+\|\nabla{\bm u}_h^{n-1}\|_{L^2}^{\frac{s}{2}} 
\|\nabla{\bm u}_h^{n}\|_{L^2}^{\frac{s}{2}} \Big) \nonumber\\
&\quad +C\tau\sum_{n=1}^N
\Big(\|{\bm u}_h^n\|_{H^1}^s
+\|{\bm f}^n\|_{H^{-1}}^s
+\|\Delta_hA_h^n\|_{L^2}^{s(5-\frac4q)}\Big) +C \nonumber \\
&\le 
C\tau\sum_{n=1}^N
\Big(
\|\nabla{\bm u}_h^{n-1}\|_{L^2}^{s(3-\frac2q)}
+\|\nabla{\bm u}_h^n\|_{L^2}^{s(3-\frac2q)}
+\|\nabla{\bm u}_h^{n-1}\|_{L^2}^s+
\|\nabla{\bm u}_h^{n}\|_{L^2}^s\Big) \nonumber\\
&\quad +C\tau\sum_{n=1}^N
\Big(\|{\bm u}_h^n\|_{H^1}^s
+\|{\bm f}^n\|_{H^{-1}}^s
+\|\Delta_hA_h^n\|_{L^2}^{s(5-\frac4q)}\Big) +C\nonumber \\
&\le 
C\tau\sum_{n=1}^N
\Big(
\|\nabla{\bm u}_h^{n-1}\|_{L^2}^2
+\|\nabla{\bm u}_h^n\|_{L^2}^2
+\|\nabla{\bm u}_h^{n-1}\|_{L^2}^2+
\|\nabla{\bm u}_h^{n}\|_{L^2}^2\Big) \nonumber\\
&\quad +C\tau\sum_{n=1}^N
\Big(\|{\bm u}_h^n\|_{H^1}^2
+\|{\bm f}^n\|_{H^{-1}}^2
+\|\Delta_hA_h^n\|_{L^2}^2\Big) +C\nonumber \\
&\le C  \qquad\forall\, s\in(1,2). 
\end{align}
Substituting \eqref{Dtau-g-1} into \eqref{Dtau-uh-g} yields 
\begin{align}\label{Dtau-u-H-1}
\tau\sum_{n=1}^N\bigg\|\frac{{\bm u}_h^n-\u_h^{n-1}}{\tau} \bigg\|_{{\bf H}^1_{0,{\rm div}}(\Omega)'}^s  
&\le C  \qquad\forall\, s\in(1,2).
\end{align}


Thus, \eqref{FEM-E}, \eqref{DtauA-L2Lp} and \eqref{Dtau-u-H-1} imply the estimates \eqref{Ahtau-uhtau} and \eqref{Ahtau-uhtau-2}. The proof of Proposition \ref{thm-energy} is complete. 
\hfill\endproof

\vspace{0.1in}
It is easy to see that the numerical scheme \eqref{numer-varphi}-\eqref{FEM-p} implies that the following equations hold 
\begin{align}
&\int_0^T\Big[\big(\mu \partial_tA_{h,\tau} ,a_h\big) 
+\big(\sigma^{-1}\nabla A_h^+,\nabla a_h\big)
- \big(\mu{\bm u}_h^+\times\big(\nabla\times A_h^-+ \mbox{$\sum_{j=1}^m$}\beta_{j,h} \nabla\times\varphi_{j,h} \big),a_h\big) \Big]\d t \nonumber \\
& 
=\int_0^T\big(\sigma^{-1} J^+ ,a_h\big) \d t   \label{FEM-A2} \\[5pt]
&\int_0^T \Big[\big(\partial_t\u_{h,\tau},{\bm v}_h\big)
+\frac12({\bm u}_h^-\cdot\nabla{\bm u}_h^+,{\bm v}_h)
-\frac12({\bm u}_h^-\cdot\nabla {\bm v}_h,{\bm u}_h^+)
+(\nu\nabla{\bm u}_h^+,\nabla {\bm v}_h)
-(p_h^+ ,\nabla\cdot{\bm v}_h)\Big]\d t \nonumber \\
&=\int_0^T
\Big[({\bm f}^+,{\bm v}_h)  + \big(\mu\big(\nabla\times A_h^-+ \mbox{$\sum_{j=1}^m$}\beta_{j,h} \nabla\times\varphi_{j,h} \big)\times \Delta_h A_h^+,{\bm v}_h\big)  \Big]\d t 
\label{FEM-u2} \\[5pt]
&\int_0^T (\nabla\cdot{\bm u}_h^+,q_h)\d t=0 \label{FEM-p2}
\end{align}
for all $a_h\in L^2(0,T;S_h^{r+1})$, 
${\bm v}_h\in L^{2}(0,T;{\bf S}_h^{r+1}) $ and $q_h\in L^2(0,T;S_h^{r})$. 

In the next subsection we pass to the limit in \eqref{FEM-A2}-\eqref{FEM-p2} by using the energy estimates given in Proposition \ref{thm-energy}.

\subsection{Compactness and convergence}\label{sec-compact}

Let $s\in(1,2)$ be fixed. 
The estimate \eqref{Ahtau-uhtau} implies the existence of functions 
\begin{align*}
&A \in L^\infty(0,T;H^1(\Omega))\cap W^{1,s}(0,T;L^2(\Omega)) \\
&
\u\in L^\infty(0,T;{\bf L}^2(\Omega)) \cap L^2(0,T;H^1(\Omega)^2)\cap W^{1,s}(0,T;H^{-1}(\Omega)^2)
\end{align*} 
such that the following convergence results hold: 
for any sequence $(h_n,\tau_n)\rightarrow 0$, there exists a subsequence $(h_{n_m},\tau_{n_m})$, $m=1,2,\dots$, such that
\begin{align}
& A_{h_{n_m},\tau_{n_m}}
&&\hspace{-10pt}\mbox{converges to $A$ weakly$^*$ in}\,\,\,L^\infty(0,T;H^1(\Omega)) , 
\label{cv-Ah}\\
& \nabla A_{h_{n_m},\tau_{n_m}}
&&\hspace{-10pt}\mbox{converges to $\nabla A$ weakly$^*$ in}\,\,\, L^\infty(0,T;{\bf L}^2(\Omega)), \label{cv-nabla-Ah}\\ 
&\partial_tA_{h_{n_m},\tau_{n_m}}
&&\hspace{-10pt}\mbox{converges to $\partial_tA$ weakly in}\,\,\, L^s(0,T;L^2(\Omega)) ,\\
& \Delta_h A_{h_{n_m},\tau_{n_m}}
&&\hspace{-10pt}\mbox{converges weakly in}\,\,\, L^2(0,T;L^2(\Omega))
\,\,\, \mbox{to a function}\,\,\, F\in L^2(0,T;L^2(\Omega)),\\ 
&\u_{h_{n_m},\tau_{n_m}}
&&\hspace{-10pt}\mbox{converges to $\u$ weakly$^*$ in}\,\,\, L^\infty(0,T;{\bf L}^2(\Omega)) 
,\\
&\u_{h_{n_m},\tau_{n_m}}
&&\hspace{-10pt}\mbox{converges to $\u$ weakly in}\,\,\, L^2(0,T;{\bf H}^1(\Omega)) 
 ,\\
&\partial_t\u_{h_{n_m},\tau_{n_m}}
&&\hspace{-10pt}\mbox{converges to $\partial_t\u$ weakly in}\,\,\, L^s(0,T;{\bf H}^1_{0,{\rm div}}(\Omega)') , \label{cv-uth}
\end{align} 
where \eqref{cv-nabla-Ah} is an obvious consequence of \eqref{cv-Ah}. Furthermore, we have $F=\Delta A \in L^2(0,T;L^2(\Omega))$. To see this, we test $F$ by a function $v\in L^2(0,T;H^1_0(\Omega))$. 
Then we have 
\begin{align*}
\int_0^T(F,v) \d t 
&= \lim_{m\rightarrow\infty} \int_0^T(\Delta_h A_{h_{n_m},\tau_{n_m}}, v) \d t \\ 
&= \lim_{m\rightarrow\infty} \int_0^T(\Delta_h A_{h_{n_m},\tau_{n_m}}, P_{h_{n_m}}v) \d t
+\lim_{m\rightarrow\infty} \int_0^T(\Delta_h A_{h_{n_m},\tau_{n_m}}, v-P_{h_{n_m}}v) \d t \\
&=-\lim_{m\rightarrow\infty} \int_0^T(\nabla A_{h_{n_m},\tau_{n_m}},\nabla P_{h_{n_m}}v) \d t
+\lim_{m\rightarrow\infty} \int_0^T(\Delta_h A_{h_{n_m},\tau_{n_m}}, v-P_{h_{n_m}}v) \d t \\
&=
 -\lim_{m\rightarrow\infty} \int_0^T(\nabla A_{h_{n_m},\tau_{n_m}},\nabla v)\d t 
 -\lim_{m\rightarrow\infty} \int_0^T(\nabla A_{h_{n_m},\tau_{n_m}},\nabla (P_{h_{n_m}}v-v)) \d t \\
&\quad+\lim_{m\rightarrow\infty} \int_0^T(\Delta_h A_{h_{n_m},\tau_{n_m}}, v-P_{h_{n_m}}v) \d t \\
&=-\int_0^T(\nabla A ,\nabla v)\d t  
 -\lim_{m\rightarrow\infty} \int_0^T(\nabla A_{h_{n_m},\tau_{n_m}},\nabla (P_{h_{n_m}}v-v)) \d t \\
&\quad+\lim_{m\rightarrow\infty} \int_0^T(\Delta_h A_{h_{n_m},\tau_{n_m}}, v-P_{h_{n_m}}v) \d t ,
\end{align*} 
where we have used \eqref{cv-nabla-Ah}. Since
\begin{align*}
\bigg|\int_0^T(\nabla A_{h_{n_m},\tau_{n_m}},\nabla (P_{h_{n_m}}v-v)) \d t \bigg|
&\le \|\nabla A_{h_{n_m},\tau_{n_m}}\|_{L^2(0,T;L^2(\Omega))}
\|\nabla (P_{h_{n_m}}v-v)\|_{L^2(0,T;L^2(\Omega))} \\
&\le C\|P_{h_{n_m}}v-v\|_{L^2(0,T;H^1(\Omega))}
\rightarrow 0\quad\mbox{as}\,\,\, m\rightarrow\infty,\\
\bigg|\int_0^T(\Delta_h A_{h_{n_m},\tau_{n_m}}, v-P_{h_{n_m}}v) \d t\bigg|
&\le \|\Delta_h A_{h_{n_m},\tau_{n_m}}\|_{L^2(0,T;L^2(\Omega))}
\|P_{h_{n_m}}v-v\|_{L^2(0,T;L^2(\Omega))} \\
&\le C\|P_{h_{n_m}}v-v\|_{L^2(0,T;L^2(\Omega))}
\rightarrow 0\quad\mbox{as}\,\,\, m\rightarrow\infty,
\end{align*} 
it follows that the function $F\in  L^2(0,T;L^2(\Omega))$ satisfies
\begin{align*}
\int_0^T(F,v) \d t 
=-\int_0^T(\nabla A ,\nabla v)\d t  \quad\forall\, v\in L^2(0,T;H^1_0(\Omega)) .
\end{align*} 
This implies $\Delta A=F\in L^2(0,T;L^2(\Omega)) $ and therefore 
\begin{align}
& \Delta_h A_{h_{n_m},\tau_{n_m}}\,\,\,\mbox{converges weakly to $ \Delta A$ in}\,\,\, L^2(0,T;L^2(\Omega))
 .
&\end{align} 

Moreover, the initial value of the limit function $A$ must be equal to $A_0$. This can be proved in the following way: on the one hand, for any given smooth function $a$ such that $a|_{t=T}=0$, 
\begin{align*}
\int_0^T(\partial_tA_{h_{n_m},\tau_{n_m}},a)\d t
&\rightarrow 
-\int_0^T(\partial_tA,a)\d t \quad\mbox{as}\,\,\, m\rightarrow\infty ; 
&\end{align*} 
on the other hand,  
\begin{align*}
\int_0^T(\partial_tA_{h_{n_m},\tau_{n_m}},a)\d t
&=-\int_0^T(A_{h_{n_m},\tau_{n_m}},\partial_ta)\d t 
+(A_{h_{n_m}}^0,a|_{t=0}) \\
&\rightarrow 
-\int_0^T(A,\partial_ta)\d t 
+(A_0,a|_{t=0}) \quad\mbox{as}\,\,\, m\rightarrow\infty.
&\end{align*} 
where we have used the convergence of $A_{h_{n_m},\tau_{n_m}}$ in \eqref{cv-Ah} and the convergence of $A_{h_{n_m}}^0$ to $A_0$ in $H^1_0(\Omega)$. The latter is a simple consequence of the convergence theory for the elliptic problem \eqref{Def-A0} and \eqref{Eq-A0h}. Therefore, 
$$-\int_0^T(\partial_tA,a)\d t =
-\int_0^T(A_{h_{n_m},\tau_{n_m}},\partial_ta)\d t 
+(A_{h_{n_m}}^0,a|_{t=0}) ,$$ 
which implies $A|_{t=0}=A_0$ in $L^2(\Omega)$. 
Similarly, we also have $\u|_{t=0}=\u_0$ in ${\bf H}^1_{0,{\rm div}}(\Omega)'$. This proves 
\begin{align}\label{ini-proved}
A|_{t=0}=A_0\quad\mbox{and}\quad \u|_{t=0}=\u_0 .
\end{align}

In Proposition \ref{thm-energy} we have proved 
$$\|\Delta_hA_{h_{n_m},\tau_{n_m}}\|_{L^2(0,T;L^2(\Omega))}+\|\partial_tA_{h_{n_m},\tau_{n_m}}\|_{L^s(0,T;L^2(\Omega))}\le C.$$
Then Lemma \ref{DeltaA-L4} and the Aubin--Lions--Simon lemma (cf. \cite[Theorem II.5.16]{BF13}) imply that $A_{h_{n_m},\tau_{n_m}}$ is compact in $L^2(0,T;W^{1,4+\alpha}(\Omega))$, and thus one can choose the subsequence to have the following property: 
\begin{align}\label{convg-A}
& A_{h_{n_m},\tau_{n_m}}\,\,\,\mbox{converges to $A$ strongly in}\,\,\,L^2(0,T;W^{1,4+\alpha}(\Omega)).
\end{align} 
Since  $A_{h_{n_m},\tau_{n_m}}-A$ is bounded in $L^\infty(0,T;H^{1}(\Omega))$ and convergent to zero in $L^2(0,T;W^{1,4+\alpha}(\Omega))$, and (interpolation inequality)
$$
\|A_{h_{n_m},\tau_{n_m}}-A\|_{L^{q_\theta}(0,T;W^{1,4}(\Omega))}
\le C\|A_{h_{n_m},\tau_{n_m}}-A\|_{L^\infty(0,T;H^{1}(\Omega))}^{1-\theta}
\|A_{h_{n_m},\tau_{n_m}}-A\|_{L^2(0,T;W^{1,4+\alpha}(\Omega))}^\theta
$$
for $\frac{1-\theta}{2}+\frac{\theta}{4+\alpha}=\frac14$ and $\frac{1}{q_\theta}=\frac{1-\theta}{\infty}+\frac{\theta}{2}$, it follows that
\begin{align}
& A_{h_{n_m},\tau_{n_m}}\,\,\,\mbox{converges to $A$ strongly in}\,\,\,L^{q_\theta}(0,T;W^{1,4}(\Omega))
\,\,\,\mbox{for some}\,\,\,q_\theta>2.
\end{align}

Since $H^1(\Omega)\hookrightarrow\hookrightarrow L^q(\Omega)\hookrightarrow {\bf H}^1_{0,{\rm div}}(\Omega)'$ for all $2\le q<\infty$ and 
$$\|\u_{h_{n_m},\tau_{n_m}}\|_{L^2(0,T;H^1(\Omega))}+\|\partial_t\u_{h_{n_m},\tau_{n_m}}\|_{L^s(0,T;{\bf H}^1_{0,{\rm div}}(\Omega)')}\le C,$$
the Aubin--Lions--Simon lemma (cf. \cite[Theorem II.5.16]{BF13}) implies that $\u_{h_{n_m},\tau_{n_m}}$ is compact in $L^2(0,T;L^q(\Omega))$ for all $2\le q<\infty$, and thus one can choose the subsequence to have the following property: 
\begin{align}\label{convg-u}
& \u_{h_{n_m},\tau_{n_m}}\,\,\,\mbox{converges to $\u$ strongly in}\,\,\,L^2(0,T;L^4(\Omega)) . 
\end{align} 

The estimates above are for the piecewise linear functions $A_{h_{n_m},\tau_{n_m}}$ and $\u_{h_{n_m},\tau_{n_m}}$. 
For the piecewise constant functions $A_{h_{n_m},\tau_{n_m}}^\pm$ and  $\u_{h_{n_m},\tau_{n_m}}^\pm$ we have similar estimates:
\begin{align}
& A_{h_{n_m},\tau_{n_m}}^\pm\,\,\,\mbox{converges to $A$ weakly$^*$ in}\,\,\,L^\infty(0,T;H^1(\Omega)) , \\
& \nabla A_{h_{n_m},\tau_{n_m}}^\pm\,\,\,\mbox{converges to $\nabla A$ weakly$^*$ in}\,\,\, L^\infty(0,T;{\bf L}^2(\Omega)), \label{AhH1-conv} \\ 
& \Delta_h A_{h_{n_m},\tau_{n_m}}^\pm\,\,\,\mbox{converges to $\Delta A$ weakly in}\,\,\, L^2(0,T;L^2(\Omega)), \label{AhDh-conv}\\ 
&\u_{h_{n_m},\tau_{n_m}}^\pm\,\,\,\mbox{converges to $\u$ weakly$^*$ in}\,\,\, L^\infty(0,T;L^2(\Omega)) ,\label{uhLinfty2-conv} \\
&\u_{h_{n_m},\tau_{n_m}}^\pm\,\,\,\mbox{converges to $\u$ weakly in}\,\,\, L^2(0,T;{\bf H}_0^1(\Omega)) 
 , \label{uhH1-conv} \\
&A_{h_{n_m},\tau_{n_m}}\,\,\,\mbox{converges to $A$ strongly in}\,\,\,L^q(0,T;W^{1,4}(\Omega))
\,\,\,\mbox{for some}\,\,\,q>2, \label{AhW14-conv} \\
&\u_{h_{n_m},\tau_{n_m}}^\pm\,\,\,\mbox{converges to $\u$ strongly in}\,\,\,L^2(0,T;L^4(\Omega)) . 
\label{uhL4-conv} 
\end{align} 
Since $\u_{h_{n_m},\tau_{n_m}}-\u$ is bounded in $L^\infty(0,T;L^2(\Omega))$ and convergent to zero in $L^2(0,T;L^4(\Omega))$, 
and (interpolation inequality)
$$
\|\u_{h_{n_m},\tau_{n_m}}-\u\|_{L^{4}(0,T;L^{\frac83}(\Omega))}
\le C\|\u_{h_{n_m},\tau_{n_m}}-\u\|_{L^\infty(0,T;L^2(\Omega))}^{1-\theta}
\|\u_{h_{n_m},\tau_{n_m}}-\u\|_{L^2(0,T;L^4(\Omega))}^\theta
$$
it follows that
\begin{align}\label{uhL83-conv}
& \u_{h_{n_m},\tau_{n_m}}\,\,\,\mbox{converges to $\u$ strongly in}\,\,\,L^4(0,T;L^{\frac83}(\Omega)) . 
\end{align} 
Then Lemma \ref{Lemma-varphij} and \eqref{uhH1-conv}-\eqref{AhW14-conv} imply
\begin{align}\label{convg-uA-nonlinear}
{\bm u}_h^+\times\big(\nabla\times A_h^-+ \mbox{$\sum_{j=1}^m$}\beta_{j,h} \nabla\times\varphi_{j,h} \big)
\,\,\,\mbox{converges to ${\bm u}\times\big(\nabla\times A + \mbox{$\sum_{j=1}^m$}\beta_j \nabla\times\varphi_{j} \big)$} \nonumber \\
\mbox{weakly in}\,\,\,L^p(0,T;L^2(\Omega)) \,\,\,\mbox{for some $p>1$} , 
\end{align}
\eqref{uhH1-conv} and \eqref{uhL83-conv} imply 
\begin{align}
{\bm u}^-\cdot\nabla{\bm u}^+\,\,\,\mbox{converges to ${\bm u}\cdot\nabla{\bm u}$ weakly in}\,\,\,L^{\frac43}(0,T;L^{\frac87}(\Omega)^2) , 
\end{align}
\eqref{uhLinfty2-conv} and \eqref{uhL4-conv} imply 
\begin{align}
{\bm u}^-\otimes{\bm u}^+\,\,\,\mbox{converges to ${\bm u}\otimes{\bm u}$ weakly in}\,\,\,L^2(0,T;L^{\frac43}(\Omega)^2) , 
\end{align} 
Lemma \ref{Lemma-varphij}, \eqref{AhDh-conv} and \eqref{AhW14-conv} imply 
\begin{align}
\big(\nabla\times A_h^-+ \mbox{$\sum_{j=1}^m$}\beta_{j,h} \nabla\times\varphi_{j,h} \big)\times \Delta_h A_h^+
\,\,\,\mbox{converges to $\big(\nabla\times A + \mbox{$\sum_{j=1}^m$}\beta_j \nabla\times\varphi_{j} \big)\times \Delta A$}\nonumber \\
\mbox{weakly in}\,\,\,L^{\frac{2q}{2+q}}(0,T;L^{\frac43}(\Omega)^2) .
\end{align}

For any given $a\in L^{\infty}(0,T;L^2(\Omega))$, its $L^2$ projection $a_h=P_ha$ converges to $a$ strongly in $L^{p'}(0,T;L^2(\Omega))$ for the number $p>1$ in \eqref{convg-uA-nonlinear}. 

Similarly, for any given ${\bm v}\in L^{\infty}(0,T;{\bf H}^1_{0,{\rm div}}(\Omega))\hookrightarrow L^{\frac43}(0,T;L^{\frac87}(\Omega)^2)'\cap L^2(0,T;L^{\frac43}(\Omega)^2)'$, its Fortin projection ${\bm v}_h=Q_h{\bm v}$ converges to ${\bm v}$ strongly in $L^{s'}(0,T;{\bf H}^1_0(\Omega))$ and satisfies $(q_h,\nabla\cdot{\bm v}_h)=0$ for all $q_h\in S_h^r$. 

Then, by taking the limit of a subsequence in \eqref{FEM-A2}-\eqref{FEM-u2}, we obtain
\begin{align}
&\int_0^T\Big[\big(\mu \partial_tA ,a\big) 
+\big(\sigma^{-1}\Delta A ,a \big)
- \big(\mu{\bm u}\times\big(\nabla\times A + \mbox{$\sum_{j=1}^m$}\beta_j \nabla\times\varphi_{j} \big),a\big) \Big]\d t  \nonumber\\
&
=\int_0^T\big(\sigma^{-1} J ,a\big) \d t  
\label{Weak-A-proved}\\
&\int_0^T \Big[\big(\partial_t\u,{\bm v}\big)
+\frac12({\bm u}\cdot\nabla{\bm u},{\bm v})
-\frac12({\bm u}\cdot\nabla {\bm v},{\bm u})
+(\nu\nabla{\bm u},\nabla {\bm v})
\nonumber   \\
&=\int_0^T 
\Big[({\bm f},{\bm v})  + \big(\mu\big(\nabla\times A + \mbox{$\sum_{j=1}^m$}\beta_j \nabla\times\varphi_{j} \big)\times \Delta A,{\bm v}\big)  \Big]\d t 
\label{Weak-u-proved}
\end{align}
for all $a\in L^{\infty}(0,T;L^2(\Omega))$ and ${\bm v}\in L^{\infty}(0,T;{\bf H}^1_{0,{\rm div}}(\Omega))$.

The convergence results \eqref{cv-Ah}-\eqref{cv-uth} show that $A$ and $\u$ possess the regularity in \eqref{RegA}-\eqref{Regu}, except the regularity ${\bm u}\in L^2(0,T;{\bf H}^1_{0,{\rm div}}(\Omega))$. 

Finally, we show that ${\bm u}\in L^2(0,T;{\bf H}^1_{0,{\rm div}}(\Omega))$. 
In fact, for any given $q\in L^2(0,T;L^2(\Omega))$, its $L^2$ projection $q_h=\widetilde P_hq$ converges to $q$ strongly in $L^2(0,T;L^2(\Omega))$. By taking limit in \eqref{FEM-p2} we obtain 
\begin{align*}
&\int_0^T (\nabla\cdot{\bm u} ,q )\d t=0 
\end{align*}
which implies $\nabla\cdot\u=0$. Since ${\bm u}\in L^2(0,T;{\bf H}^1_0(\Omega))$ and $\nabla\cdot\u=0$, it follows that ${\bm u}\in L^2(0,T;{\bf H}^1_{0,{\rm div}}(\Omega))$. 

Therefore, \eqref{ini-proved} and \eqref{Weak-A-proved}-\eqref{Weak-u-proved} imply that $(A,\u,(\varphi_j)_{j=1}^m,(\beta_j)_{j=1}^m)$ is a weak solution of \eqref{PDE-A}-\eqref{Def-A0}. 
Furthermore, \eqref{convg-A}, \eqref{convg-u} and Lemma \ref{Lemma-varphij} imply the convergence results in Theorem \ref{MainTHM}.  
\hfill\endproof

\section{Numerical examples}

In this section we present two numerical examples to illustrate the convergence of the proposed numerical method in nonconvex and nonsmooth domains. All the numerical tests are done by using FreeFEM++ with double precision. 

\begin{example}\label{Example:Lshape}
{\upshape In the first example, we consider the MHD equations 
\begin{align} 
&\partial_t {\bm H}
+\nabla\times(\nabla\times{\bm H})
-\nabla\times({\bm u}\times{\bm H})
=\nabla\times J  \label{Exp-eq1} \\ 
&\partial_t{\bm u} +{\bm u}\cdot\nabla{\bm u}-\Delta{\bm u}
+\nabla p 
={\bm f} - \mu {\bm H}\times(\nabla\times{\bm H})  \label{Exp-eq2} \\
&\nabla\cdot{\bm u}=g \label{Exp-eq3}  
\end{align} 
in a simply connected L-shape domain $\Omega$ whose longest side has unit length, centered at the origin; see Figure \ref{Lshape} (left). 
The source terms ${\bm f}$, $g$ and 
$$
J= \partial_tu_1+ \nabla\times{\bm H}-{\bm u}\times{\bm H} 
$$ 
are calculated by substituting the following exact solution into \eqref{Exp-eq1}-\eqref{Exp-eq3}: 
\begin{align}
&\u=\left(\begin{aligned}
u_1\\
u_2
\end{aligned}\right)
=\left(\begin{aligned}
t^2\Phi(r)r^{2/3}\sin(2\theta/3)\\
t^2\Phi(r)r^{2/3}\sin(2\theta/3)
\end{aligned}\right) ,
\quad p=0, \label{sol-u} \\
&{\bm H}=
\nabla\times u_1=
\left(\begin{aligned} 
2t^2 \Phi(r)r^{-1/3}/3\cos(\theta/3)
+t^2\Phi'(r)r^{2/3} \sin(2\theta/3)\sin(\theta) 
 \\
2t^2 \Phi(r)r^{-1/3}/3 \sin(\theta/3 )
-t^2\Phi'(r)r^{2/3} \sin(2\theta/3)\cos(\theta)
\end{aligned}\right) , \label{sol-H} 
\end{align}
where {\color{blue}$r=\sqrt{x^2+y^2}$ and $\theta={\rm arg}(x+iy)$, with $(x,y)\in\Omega$,} 
and $\Phi(r)$ is a $C^3(\Omega)$ cut-off function. Here, 
$$
\Phi(r)=\left\{
\begin{array}{ll}
0.1 & \mbox{if}~~r<0.1  \\
\Upsilon(r) &\mbox{if}~~ 0.1\leq r\leq 0.4 \\
0 & \mbox{if}~~r>0.4
\end{array}\right.
$$
and $\Upsilon(r)$ is the unique $7^{\rm th}$ 
order polynomial satisfying the 
conditions $\Upsilon'(0.1)=\Upsilon''(0.1)
=\Upsilon'''(0.1)=\Upsilon(0.4)=\Upsilon'(0.4)
=\Upsilon''(0.4)=\Upsilon'''(0.4)=0$ 
and $\Upsilon(0.1)=0.1$. The constructed solutions in \eqref{sol-u}-\eqref{sol-H} satisfy the boundary conditions 
\begin{align} 
\H\cdot{\bm n}=0,
\quad 
\nabla\times\H = J ,
\quad\mbox{and}\quad
\u=0 
\quad\mbox{on}\,\,\,\Omega\times(0,T] . 
\end{align} 

We solve the equations \eqref{Exp-eq1}-\eqref{Exp-eq3} up to time $T=1$ by the proposed numerical method \eqref{numer-varphi}-\eqref{FEM-p} with $r=1$, i.e., with the P2 element for $A$ and P2-P1 elements for $({\bm u},p)$. 
The numerical solution of magnetic field is given by 
$
{\bm H}_h^n=\nabla\times A_h^{n} . 
$
Since the L-shape domain is simply connected, it follows that $m=0$ (the constants $\beta_j$, $j=1,\dots,m$, are not needed). The L-shape domain is triangulated quasi-uniformly with $M$ nodes per unit length on each side, and we denote by $h=1/M$ for simplicity. 

\begin{figure}[t]
\begin{center}
\begin{tikzpicture}
\draw[fill=gray!15] (0,0) -- (1.5,0) -- (1.5,1.5) -- (-1.5,1.5) -- (-1.5,-1.5) -- (0,-1.5) -- (0,0) ; 
\draw (1.5,0) node[below]{\footnotesize$(0.5,0)$} ;  
\draw (1.5,1.5) node[above]{\footnotesize$(0.5,0.5)$} ;  
\draw (-1.5,1.5) node[above]{\footnotesize$(-0.5,0.5)$} ;  
\draw (-1.5,-1.5) node[below]{\footnotesize$\,\,\,\,(-0.5,-0.5)$} ;  
\end{tikzpicture}
\hspace{30pt}
\begin{tikzpicture}
\draw[fill=gray!15] (-1.5,-1.5) -- (1.5,-1.5) -- (1.5,1.5) -- (-1.5,1.5) -- (-1.5,-1.5) ; \draw[fill=white] (0.5,-0.5) -- (0.5,0.5) -- (-0.5,0.5) -- (-0.5,-0.5) -- (0.5,-0.5); 
\draw (-1.5,1.5) node[above]{\footnotesize$(-0.5,0.5)$} ;  
\draw (-1.5,-1.5) node[below]{\footnotesize$(-0.5,-1)$} ;  
\draw (1.5,-1.5) node[below]{\footnotesize$(1,-1)$} ;  
\draw (-0.5,0.5) node[above]{\hspace{-3pt}\footnotesize$(0,0)$} ;  
\draw (-0.5,-0.5) node[below]{\footnotesize$(0,-0.5)$} ;  
\draw (0.5,0.5) node[above]{\hspace{3pt}\footnotesize$(0.5,0)$} ;  
\end{tikzpicture}
\end{center}
\caption{\small An L-shape domain (Example \ref{Example:Lshape}) and a multi-connected domain (Example \ref{Expl-2})}
\label{Lshape}
\end{figure}
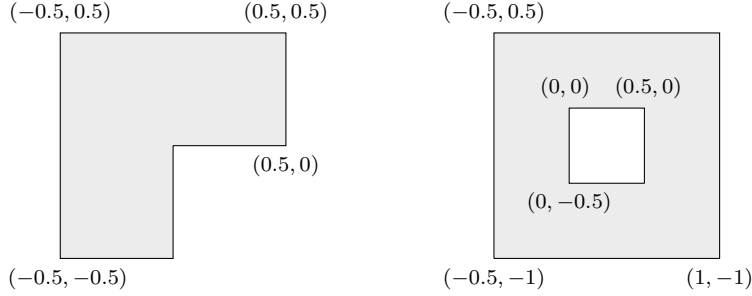

We compare the numerical solutions with the exact solution given by \eqref{sol-u}-\eqref{sol-H} and present the errors of the numerical solutions in Table  \ref{Example1-1}. For comparison, we also present the numerical results of a ``direct $H^1$-conforming FEM'' in Table \ref{Example1-2}, with the same time-stepping scheme as the proposed method \eqref{numer-varphi}-\eqref{FEM-p}. In particular, the direct $H^1$-conforming FEM seeks ${\bm H}_h^n\in {\bf S}_h^2$, $\u_h^n\in \mathring {\bf S}_h^2$ 
and $p_h^n\in S_h^1$, with ${\bm H}_h^n\cdot{\bm n}=0$ on $\partial\Omega$, such that the following equations hold for all test functions ${\bf a}_h\in {\bf S}_h^2$, ${\bm v}_h\in \mathring{\bf S}_h^2$ and $q_h\in S_h^1$ with ${\bf a}_h\cdot{\bm n}=0$ on $\partial\Omega$:
\begin{align}\label{H1-conf-1}
&\bigg(\mu\frac{{\bm H}_h^n-{\bm H}_h^{n-1}}{\tau} ,{\bf a}_h\bigg)
+(\sigma^{-1}\nabla \times {\bm H}_h^n,\nabla \times {\bf a}_h)
- (\mu{\bm u}_h^{n}\times {\bm H}_h^{n-1} , \nabla\times {\bf a}_h) \nonumber \\
&
=(\sigma^{-1} J^n ,\nabla\times {\bf a}_h)  \\[8pt]
&\bigg(\frac{{\bm u}_h^n-\u_h^{n-1}}{\tau} ,{\bm v}_h\bigg)
+\frac12({\bm u}_h^{n-1}\cdot\nabla{\bm u}_h^n,{\bm v}_h)
-\frac12({\bm u}_h^{n-1}\cdot\nabla {\bm v}_h,{\bm u}_h^n)
+(\nu\nabla{\bm u}_h^n,\nabla {\bm v}_h)
-(p_h^n ,\nabla\cdot{\bm v}_h) \nonumber \\
&=({\bm f}^n,{\bm v}_h)  - (\mu{\bm H}_h^{n-1} \times (\nabla\times{\bm H}_h^n),{\bm v}_h)   
 \\[8pt]
&(\nabla\cdot{\bm u}_h^n,q_h)=(g^n,q_h) .  \label{H1-conf-3}
\end{align}

The numerical results in Tables \ref{Example1-1} and \ref{Example1-2} show that the proposed method has slightly higher accuracy than the  $H^1$-conforming FEM in computing the magnetic field ${\bm H}$ (with the same degree of finite elements and similar computational complexity). The convergence of the proposed numerical method is proved in Theorem \ref{MainTHM} (also see Remark \ref{Remark-THM}), while the convergence of the direct $H^1$-conforming FEM remains open in nonconvex and nonsmooth domains.

\begin{table}[!htbp]
\centering
\caption{\small Numerical results of the proposed numerical method \eqref{numer-varphi}-\eqref{FEM-p} 
with $\tau=h$} \vspace{-5pt}
\begin{center}
\begin{tabular}{c|c|cc}\hline
$h$ 
& $\|{\bm H}_h^N-{\bm H}(\cdot,t_N)\|_{L^2}$ & $\|{\bm u}_h^N-{\bm u}(\cdot,t_N)\|_{L^2}$ \\ \hline
1/16   &4.295E-03  &6.540E-05  \\ \hline
1/32   &1.894E-03  &2.073E-05 \\ \hline   
1/64   &1.092E-03  &8.492E-06 \\ \hline   
1/128 &7.036E-04  &3.845E-06 \\ \hline   
convergence rate
&$O(h^{0.63})$ &$O(h^{1.14})$ \\ \hline
\end{tabular}
\end{center}
\label{Example1-1}
\end{table}

\vspace{-5pt}
\begin{table}[!htbp]
\centering
\caption{\small Numerical results of the $H^1$-conforming FEM \eqref{H1-conf-1}-\eqref{H1-conf-3} with $\tau=h$}\vspace{-5pt}
\begin{center}
\begin{tabular}{c|c|c}\hline
$h$ & $\|{\bm H}_h^N-{\bm H}(\cdot,t_N)\|_{L^2}$ & $\|{\bm u}_h^N-{\bm u}(\cdot,t_N)\|_{L^2}$ \\ \hline
1/16   &1.180E-02  &6.516E-05   \\ \hline 
1/32   &7.419E-03  &2.046E-05   \\ \hline   
1/64   &4.958E-03  &8.297E-06   \\ \hline   
1/128 &3.206E-03  &3.723E-06    \\ \hline  
convergence rate
&$O(h^{0.63})$  &$O(h^{1.16})$  \\ \hline
\end{tabular}
\end{center}
\label{Example1-2}
\end{table}

}
\end{example}

\begin{example}\label{Expl-2}
{\upshape 
In the second example, we consider the MHD equations \eqref{Exp-eq1}-\eqref{Exp-eq3} in a multi-connected domain shown in Figure \ref{Lshape} (right), with the source terms ${\bm f}$, $g$ and 
$$
J= \partial_tu_1+ \nabla\times{\bm H}-{\bm u}\times{\bm H} 
$$ 
determined by the exact solution 
\begin{align}
&\u=\left(\begin{aligned}
u_1\\
u_2
\end{aligned}\right)
=\left(\begin{aligned}
t^2\Phi(r)r^{2/3}\sin(2\theta/3)\\
t^2\Phi(r)r^{2/3}\sin(2\theta/3)
\end{aligned}\right) ,
\quad p=0, \label{sol-u-2} \\
&{\bm H}=
\nabla\times u_1+\nabla\times\varphi , \label{sol-H-2} 
\end{align}
where $\varphi$ is the harmonic function satisfying 
\begin{align}\label{PDE-varphi}
\left\{\begin{aligned}
&\Delta \varphi=0 &&\mbox{in}\,\,\,\Omega \\
&\varphi=1&&\mbox{on the inner boundary} \\
&\varphi=0 &&\mbox{on the outer boundary}. 
\end{aligned}\right.
\end{align}
In this case $m=1$ and $\beta_1=1$.

We solve the MHD equations up to time $T=1$ by the proposed numerical method \eqref{numer-varphi}-\eqref{FEM-p} with the P2 element for $A$ and P2-P1 elements for $({\bm u},p)$, and compare the numerical solutions with the exact solution given by \eqref{sol-u-2}-\eqref{sol-H-2} (where $\varphi$ is approximated numerically). The domain is triangulated quasi-uniformly, with $M$ nodes per unit length on each side, and we denote by $h=1/M$ for simplicity. 
The errors of the numerical solutions are presented in Table  \ref{Example2-1}. For comparison, we also present the numerical results of the direct $H^1$-conforming FEM \eqref{H1-conf-1}-\eqref{H1-conf-3} with the P2 element for ${\bm H}$ and P2-P1 elements for $({\bm u},p)$ in Table \ref{Example2-2}. 
Numerical results in Tables \ref{Example2-1}--\ref{Example2-2} show that the proposed method is much more accurate than the direct $H^1$-conforming FEM in such a multi-connected nonsmooth domain. The reason may be that the proposed method approximates the harmonic part $\nabla\times\varphi$ at the initial time in a more accurate way than the direct $H^1$-conforming FEM, which cannot separate $\nabla\times\varphi$ from ${\bm H}$. 



\begin{table}[!htbp]
\centering
\caption{\small Numerical results of the proposed numerical method \eqref{numer-varphi}-\eqref{FEM-p} 
with $\tau=h$}\vspace{-5pt}
\begin{center}
\begin{tabular}{c|c|cc}\hline
$h$ 
& $\|{\bm H}_h^N-{\bm H}(\cdot,t_N)\|_{L^2}$ & $\|{\bm u}_h^N-{\bm u}(\cdot,t_N)\|_{L^2}$ \\ \hline
1/16  &5.625E-02  &6.376E-05 \\ \hline
1/32  &3.472E-02  &2.232E-05 \\ \hline   
1/64  &2.179E-02  &8.678E-06 \\ \hline   
1/128 &1.366E-02 &3.767E-06 \\ \hline  
convergence rate
&$O(h^{0.67})$ &$O(h^{1.20})$ \\ \hline
\end{tabular}
\end{center}
\label{Example2-1}
\end{table}

\vspace{-5pt}
\begin{table}[!htbp]
\centering
\caption{\small Numerical results of the $H^1$-conforming FEM \eqref{H1-conf-1}-\eqref{H1-conf-3} with $\tau=h$}\vspace{-5pt}
\begin{center}
\begin{tabular}{c|c|c}\hline
$h$ & $\|{\bm H}_h^N-{\bm H}(\cdot,t_N)\|_{L^2}$ & $\|{\bm u}_h^N-{\bm u}(\cdot,t_N)\|_{L^2}$ \\ \hline
1/16   &4.681E-01  &1.290E-04   \\ \hline   
1/32   &2.912E-01  &4.582E-05   \\ \hline   
1/64   &1.819E-01  &1.798E-05   \\ \hline   
1/128 &1.153E-01  &7.388E-06    \\ \hline  
convergence rate
&$O(h^{0.66})$  &$O(h^{1.28})$  \\ \hline
\end{tabular}
\end{center}
\label{Example2-2}
\end{table}

}
\end{example}

{\color{blue}
\section{Conclusion and remarks}

We have proposed a fully discrete and linearized $H^1$-conforming Lagrange finite element method for solving a magnetic potential formulation of the two-dimensional MHD equations in general domains that can be non-convex, nonsmooth and multi-connected.  The theoretical analysis shows that, for given source terms and initial data with regularity \eqref{Reg-f-J}, a subsequence of numerical solutions converges in $L^2(0,T;{\bf L}^2(\Omega))$ to a weak solution of the MHD equations without any mesh restriction.  
The uniqueness of weak solutions for the two-dimensional MHD equations in general domains remains open. If the weak solution is unique, then the numerical solutions converge to the unique weak solution as $(\tau,h)\rightarrow (0,0)$ (without taking a subsequence). 

The results proved in this paper can be generalized to more general $H^1$-conforming finite element spaces. For example, the finite element space $\mathring {\bf S}_h^{r+1}\times S_h^r/\R$ for $(\u,p)$ can be replaced by other $H^1$-conforming finite element spaces $X_h\times M_h$ satisfying the inf-sup condition (4.1), meanwhile with the following approximation property: 
$$
\lim_{h\rightarrow 0}
\inf_{(\v_h,q_h)\in X_h\times M_h}
(\|\v-\v_h\|_{H^1}+\|q-q_h\|_{L^2})
=0,\quad \forall\, (\v,q)\in \H^1_0(\Omega)\cap L^2(\Omega).
$$
The finite element space $\mathring S_h^{r+1}$ for $A$ can be replaced by any $H^1$-conforming finite element space $V_h\subset H^1_0(\Omega)$ such that 
$$
\inf_{\chi_h\in V_h}\|\chi-\chi_h\|_{H^1}\le Ch^\alpha\|\chi\|_{H^{1+\alpha}},\quad \forall\, \chi\in H^1_0(\Omega)\cap H^{1+\alpha}(\Omega),\,\,\,\forall\,\alpha\in[0,1].
$$

The assumption of perfect-conductor-type boundary conditions for the magnetic field is an ideal situation. In realistic problems, the magnetic field is often not confined to the fluid region, but extends throughout space, satisfying interface conditions (rather than boundary conditions) on the surface of the fluid region. In this case, the numerical approximation of solutions to the full incompressible MHD equations in the whole space, satisfying interface conditions on the surface of the fluid region, is still an open problem. The present paper, focused on two-dimensional domains with perfect-conductor-type boundary conditions, could be an incremental step towards the solution. 
}

\bigskip

{\bf Funding.} The work of Buyang Li was supported in part by the Hong Kong RGC grant 15301818. The work of Jilu Wang was supported in part by the USA National Science Foundation grant DMS-1315259 and by the USA Air Force Office of Scientific Research grant FA9550-15-1-0001. The work of Liwei Xu was supported in part by a Key Project of the Major Research Plan of NSFC grant no.~91630205 and the NSFC grant no.~11771068. 


\bibliographystyle{amsplain}

\end{document}